\documentclass{elsarticle}

\usepackage{hyperref}
\usepackage{bm}
\usepackage{float}
\usepackage{amsmath}
\usepackage{amsthm}
\usepackage{amsfonts}
\usepackage{graphicx}
\usepackage{comment}

\usepackage[usenames,dvipsnames]{xcolor}

\newtheorem{algorithm}{Algorithm}[section]

\newtheorem{remark}{Remark}[section]

\usepackage{color}

\newcommand{\rev}[1]{\textcolor{black}{{{#1}}}}
\usepackage{setspace}
\begin{document}

\begin{frontmatter}

\title{Local time stepping for the shallow water equations in \rev{MPAS}}

\author{Giacomo Capodaglio \fnref{lanlFootnote}}

\author{Mark Petersen \fnref{lanlFootnote}}

\fntext[lanlFootnote]{Computational Physics and Methods Group
Los Alamos National Laboratory,
e-mails: gcapodaglio, mpetersen@lanl.gov}

\begin{abstract}
\rev{We assess the performance of a set of local time-stepping (LTS) schemes for the shallow water equations implemented in the Model for Prediction Across Scales (MPAS). The goal of LTS is to speed up the simulation by allowing different time-steps on different regions of the computational grid. The LTS schemes considered here were originally introduced by Hoang et al. (J. Comput. Phys., Vol. 382, p.152-176, 2019), who laid out the mathematical foundation of the methods. Here, the authors take on the task of presenting a fast, efficient and scalable parallel implementation of these LTS methods on high performance computing machines, with the aim to provide a recipe for other climate modeling groups that may be interested in employing LTS algorithms in their codes. 
As a matter of fact, even if MPAS is our framework of choice, our approach is general enough and could be of interest to other groups beyond the MPAS community.
Due to their nature, LTS methods possess an inherent load imbalance that needs to be carefully addressed in order to obtain efficient scalability. Even more important is the far from trivial task of computing the right-hand side terms only on specific LTS regions during the time-stepping procedure. An inefficient handling of this task causes a drastic decay of the CPU time performance, making the LTS algorithms practically of no use.
The emphasis of the present work is therefore on the computational and parallel aspects of the LTS methods, whose proper handling is crucial to make the methods run faster against existing strategies, such as for instance high-order explicit global time-stepping schemes. This is in fact the ultimate goal of using an LTS procedure and it is the one to which we direct all our optimization efforts.}
\end{abstract}
 

\end{frontmatter}

\section{Introduction}
\rev{The Model for Prediction Across Scales (MPAS) is a climate model framework developed by Los Alamos National Laboratory (LANL) in collaboration with the National Center for Atmospheric Research (NCAR), which includes various components such as ocean \cite{ringler2013multi},  sea-ice \cite{turner2018mpas}, land-ice \cite{hoffman2018mpas}, and atmosphere \cite{skamarock2012multiscale}. The MPAS ocean and ice components are part of the Energy Exascale Earth System Model (E3SM), developed by the U.S. Department of Energy, which runs full climate simulations on variable-resolution meshes \cite{Golaz_2019_JAMES, Petersen_2019_JAMES, caldwell2019}.  A shallow water core is included in MPAS \cite{ringler2011momentum} for the development and testing of fundamental algorithms, and is the focus of the present work.
Within the MPAS framework, the horizontal discretization for the shallow water core is characterized by a C-grid staggering finite-volume type scheme referred to as the TRiSK scheme \cite{thuburn2009numerical, ringler2010unified}, which is up to second order accurate in space. The TRiSK scheme is built on a spherical centroidal Voronoi tessellation (SCVT) mesh \cite{du2002grid,du1999centroidal}, which allows for user-defined variable resolution through the definition of appropriate mesh-density functions that dictate the size of the grid within the domain \cite{ju2011voronoi}. The Voronoi tessellation is referred to as the primal mesh, whose dual is an associated Delaunay triangulation, where the vertices coincide with the centers of the Voronoi cells. Both the primal and the dual grids are available in MPAS, and they are currently built using the JIGSAW library \cite{engwirda2017jigsaw};
see \cite{hoch2020mpas} for a recent study on the capabilities of JIGSAW in MPAS-Ocean. 
Although the vertical discretization will not be immediately relevant to the work presented in this paper, we refer to \cite{petersen2015evaluation} for an in-depth analysis of vertical coordinate types within MPAS-Ocean.}

\rev{The option to perform simulations on variable resolution SCVT type meshes is a distinctive trait of MPAS that makes it different from existing models. A detailed study on the features of a multi-resolution ocean model has been presented in \cite{ringler2013multi}. In a multi-resolution context, the choice of time-step for explicit time integrators has to respect an upper bound that depends on the size of the smallest cell on the grid, according to the Courant-Friedrichs-Lewy (CFL) condition. This represents a reason for concern because the small time-step of the high resolution region would have to be used also on the low resolution region, where the physical processes could allow much larger time-steps, resulting in an unnecessary computational effort and longer simulation times.}

To overcome issues of this kind, local time-stepping (LTS) schemes have been developed, with the aim of allowing a different choice of time-steps in different regions of the mesh, so that larger time-steps could be used where small time-steps are not required.
Several LTS schemes are available in the literature, and providing a complete overview is beyond the scope of this paper. 
To name a few, in \cite{maleki2016novel, trahan2012local} the authors consider a discontinuous Galerkin type of discretization for the shallow water equations, while in \cite{dazzi2016local, fumeaux2004generalized} a finite volume approach is discussed.
Other works have addressed local time-stepping for conservation laws \cite{osher1983numerical, bremer2020adaptive,tan2004moving}, wave propagation \cite{grote2015runge} and mixed flows \cite{dazzi2016local}. For an example of work that tackles specifically coastal applications we refer to \cite{dawson2013parallel}.
The LTS work that is most relevant to the MPAS framework is perhaps \cite{hoang2019conservative}, where the authors developed second and third order schemes based on the strong stability preserving Runge-Kutta (SSPRK) methods for the shallow water equations discretized with TRiSK on SCVT grids. 
\rev{ The authors of \cite{hoang2019conservative} presented the mathematical foundation of their LTS algorithms, showing desirable conservation properties of the methods as well as appropriate numerical behavior. The goal of the present work is to provide a detailed discussion on how to make the right design choices to efficiently implement in parallel the LTS schemes in \cite{hoang2019conservative}, in a way that would make them faster than high-order explicit global time-stepping schemes such as for instance the well-known Runge-Kutta scheme of fourth order. How to achieve this goal is far from trivial and was not discussed in \cite{hoang2019conservative}.
Among the challenges that have to be addressed, there is the inherent load imbalance of the LTS schemes due to the fact that different time-steps are used on different regions of the mesh. Moreover, achieving good parallel scalability does not by itself guarantee that the LTS schemes will be faster than explicit global time-stepping schemes, a feature that is the most important proof of the usefulness of the LTS strategy, and a necessary step to convince the scientific community that LTS is indeed a valuable tool for multi-resolution grids. As a matter of fact, one of the biggest challenges is to find a way to efficiently handle the parallel partitioning at the same time as the LTS partitioning, in order to compute the right-hand side terms only on the specific LTS region in which the solution should be advanced during the time-stepping procedure.
We will expand on these topics in the remainder of the paper.
We observe that even though we are using MPAS as our framework of choice, the approach is general enough and could be of interest to other groups beyond the MPAS community.}

\rev{The paper is structured as follows: 
in Section \ref{LTS} we report the LTS schemes from \cite{hoang2019conservative}, as well as the SSPRK algorithms on which they are based. The LTS schemes have been recast in a different shape in our implementation in MPAS to minimize storage requirements and address some technicalities that were not fully considered in \cite{hoang2019conservative}.
In Section \ref{numRes} we present a detailed computational analysis of the performance of the LTS algorithms described in the previous section, with a particular focus on parallel performance, and CPU times.
Finally, in Section \ref{End} we draw our conclusions and discuss future work.
For the reader's references, we included in Appendix A a brief review of the nonlinear shallow water equations, whereas their spatial discretization using the TRiSK scheme on SCVT grids is reported in Appendix B.}
%%%%%%%%%%%%%%%%%%%%%%%%%%%%%%%%%%%%%%%%%%%%%%%%

\section{Time stepping schemes}\label{LTS}
 
\rev{Let us consider the discretized nonlinear shallow water equations:}
\rev{
\begin{equation}
\begin{cases}
\begin{aligned}\label{eq:sweTimeStepping}
    \dfrac{\partial h_i}{\partial t} &= - \Big[\nabla \cdot F_e\Big]_i := \mathcal{H}_i(\mathbf{h},\mathbf{u}),\\
    \dfrac{\partial u_e}{\partial t} &= - F_e^{\perp}[q]_{v \rightarrow e}  - \Big[g \nabla(h_i +b_i) +\nabla K_i\Big]_e:= \mathcal{U}_e(\mathbf{h},\mathbf{u}),
\end{aligned}
\end{cases}
\end{equation}
}
where $\mathbf{h}=\{h_i\}$ and $\mathbf{u}=\{u_e\}$ are the vectors of center and edge values, respectively.
\rev{As a reference, the formulation of the shallow water equations system as well as its discretization with the TRiSK scheme can be found respectively in Appendix A and Appendix B at the end of this document, together with definitions of symbols and variables.}
Let $T>0$ and partition the time interval $[0,T]$ in $N$ sub-intervals of size $\Delta t = T/N$ such that $[0,T]=\cup_{n=1}^N[t_{n-1},t_n]$, with $t_0=0$ and $t_n=t_{n-1}+\Delta t$, for $n=1,\ldots N$.
Assume that the values of $(\mathbf{h},\mathbf{u})$ at $t_0$ are known, denote them by $(\mathbf{h}^0,\mathbf{u}^0)$.  From now on, for a given time level $t_n$, let $(\mathbf{h}^{old},\mathbf{u}^{old})$ be the vectors at time level $t_{n}$ and $(\mathbf{h}^{new},\mathbf{u}^{new})$ the vector at the next time level, $t_{n+1}$. 

\subsection{Strong stability preserving Runge-Kutta methods}
The LTS schemes in \cite{hoang2019conservative} are based on explicit methods of second and third order accuracy, normally referred to as strong stability preserving Runge-Kutta (SSPRK) methods \cite{gottlieb2001strong, shu1988total}.
The second order SSPRK method, from now on denoted SSPRK2, is a two-stage method, whereas the third order SSPRK method, from now on denoted SSPRK3, is a three-stage method. The two methods can be written compactly as one with the following notation.

\begin{algorithm}[SSPRK]
Assume $(\mathbf{h}^{old},\mathbf{u}^{old})$ and $\Delta t$ are given. Then, depending on whether SSPRK2 or SSPRK3 is used, $(\mathbf{h}^{new},\mathbf{u}^{new})$ is computed with the following stages:

\underline{Stage 1}.
\begin{equation}
\begin{aligned}
\begin{cases}
h_i^{1st} = h_i^{old} + \Delta t \, \mathcal{H}_i(\mathbf{h}^{old},\mathbf{u}^{old})\\
u_e^{1st} = u_e^{old} + \Delta t \, \mathcal{U}_e(\mathbf{h}^{old},\mathbf{u}^{old})
\end{cases}.
\end{aligned}
\end{equation}
\underline{Stage 2}.
\begin{equation}
\begin{aligned}
\begin{cases}
h_i^{2nd} = \omega^{old} h_i^{old} + \omega^{1st} h_i^{1st} + \omega^{rhs} \Delta t \, \mathcal{H}_i(\mathbf{h}^{1st},\mathbf{u}^{1st})\\
u_e^{2nd} = \omega^{old} u_e^{old} + \omega^{1st} u_e^{1st} + \omega^{rhs} \Delta t \, \mathcal{U}_e(\mathbf{h}^{1st},\mathbf{u}^{1st})\\
\end{cases}.
\end{aligned}
\end{equation}
\begin{equation}
\omega^{old}=
    \begin{cases}
    1/2, \quad \mbox{if SSPRK2}\\
    3/4, \quad \mbox{if SSPRK3}
    \end{cases},\qquad 
    \omega^{1st}=\omega^{rhs}=
    \begin{cases}
        1/2, \quad \mbox{if SSPRK2}\\
    1/4, \quad \mbox{if SSPRK3}
    \end{cases}.
\end{equation}
If SSPRK2: 
\begin{equation}
\begin{aligned}
\begin{cases}
h_i^{new} = h_i^{2nd}  \\
u_e^{new} = u_e^{2nd}  
\end{cases}.
\end{aligned}
\end{equation}
Else if SSPRK3: \\
\underline{Stage 3}.
\begin{equation}
\begin{aligned}
\begin{cases}
h_i^{new} = \dfrac{1}{3}h_i^{old} + \dfrac{2}{3}h_i^{2nd} + \dfrac{2}{3}\Delta t \,\mathcal{H}_i(\mathbf{h}^{2nd},\mathbf{u}^{2nd})   \\
u_e^{new} = \dfrac{1}{3}u_e^{old} + \dfrac{2}{3}u_e^{2nd} + \dfrac{2}{3}\Delta t \,\mathcal{U}_e(\mathbf{h}^{2nd},\mathbf{u}^{2nd})
\end{cases}.
\end{aligned}
\end{equation}
End if.
\end{algorithm}

As remarked in \cite{hoang2019conservative}, the  SSPRK3 method permits larger time-steps compared to SSPRK2, and it is more efficient and practically more stable.

\subsection{Local time stepping schemes}
\begin{figure}[!t]
   \centering
   \includegraphics[scale=0.6]{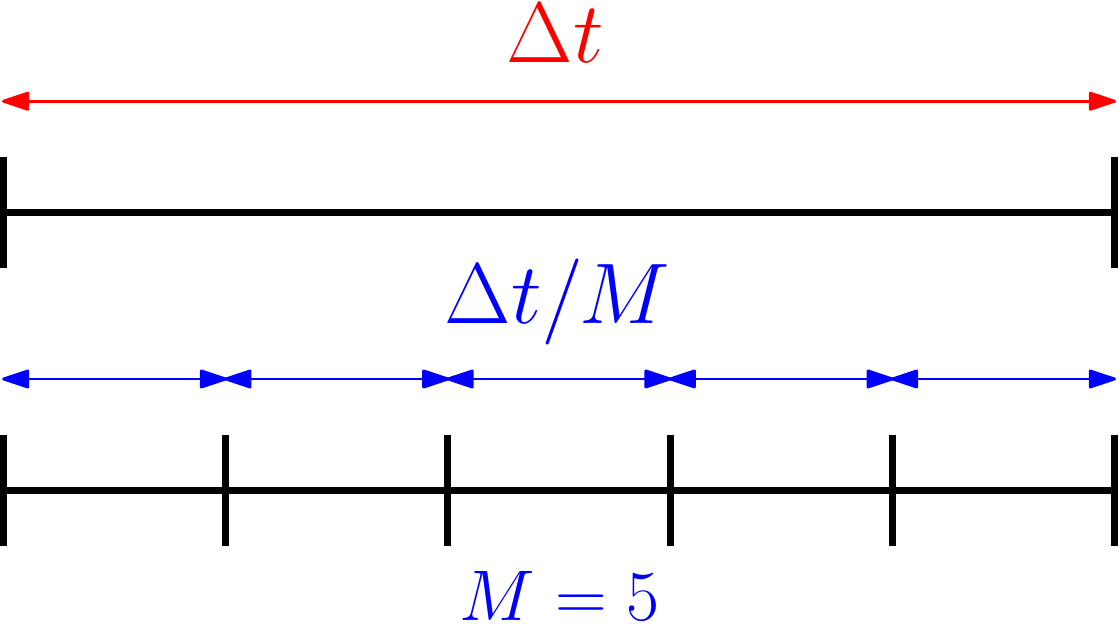}
\caption{Partitioning of the time interval according to coarse and fine time increments.}
   \label{fig:2}
\end{figure}

The LTS schemes in \cite{hoang2019conservative} have been built with the intention of allowing the selection of different time-step sizes depending on the specific resolution of different regions of the mesh.
Ideally, the partitioning of the mesh would be done so that the mesh is divided in two regions, a coarse time-step region, where the numerical integration is carried out using the coarse time-step $\Delta t = T/N$, and a fine time-step region, where a fine time-step $\Delta t / M = T / (NM)$ is used, with $M \in \mathbb{N}^+$, see Figure \ref{fig:2} for a sketch.
%%%%%%%%%%%%
%%%%%%%%%%%%
Then, to handle the communication between the two regions, an interface layer is placed between them.
\rev{For the LTS schemes developed in \cite{hoang2019conservative}, 
there are two interface layers, the second of which we will refer to as interface layer 2, needed exclusively to guarantee exact mass conservation, see Theorem 4.4 in \cite{hoang2019conservative}.}
Moreover, due to the way the algorithms have been designed, the coarse time-step $\Delta t$ is supposed to be used also within the two interface regions.
Technically, a decision on how to assign a specific LTS region to a cell of the mesh should be made according to considerations on the cell size, in order to decide what time-step to use on that specific cell based on a CFL condition.
However, we feel it is important to remember that the goal of LTS is to use large time-steps on {\it most} of the low resolution region, in order not to slow down the computational time. Hence, because the LTS schemes in \cite{hoang2019conservative} are designed so that one has to use the coarse time-step on the interface regions, it may be helpful to select these regions not necessarily according to cell dimensions but rather on making sure that the size of their cells does not pose a constraint on the maximum size of the time-step that can be used sufficiently far away from the transition region. 
Specifically, it may be beneficial to make up the interface layers with cells that have approximately the same size as the coarse region cells, so that it is certain that the maximum coarse time step can be selected. This comment provides practical insight on the use of the LTS methods in \cite{hoang2019conservative} and exemplifies how our work is different from  \cite{hoang2019conservative}, which was intended to show the mathematical formulation and conservation properties of the methods, without necessarily digging too deep on their implementation, practical aspects and parallel performance. Exploring these areas is the goal of the present paper.

\begin{figure}[!t]
   \centering
   \includegraphics[scale=1.8]{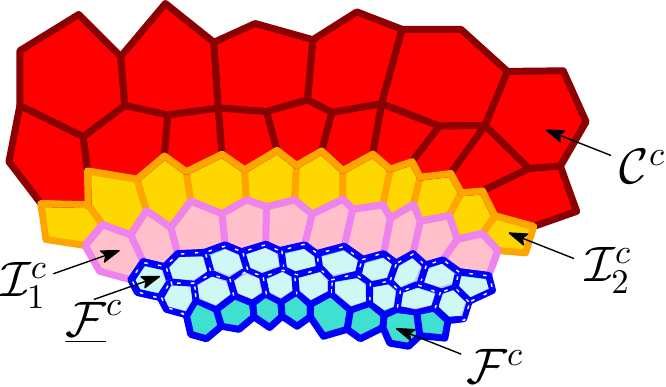}
\caption{LTS regions: fine (blue), interface 1 (pink), interface 2 (yellow) and coarse (red). For ease of visualization, the edges in $\underline{\mathcal{F}^e}$ are described by a white dotted line.}
   \label{fig:3}
\end{figure}

For the sake of laying out the LTS algorithms,
let us consider the following mutually disjoint index sets that are associated with cell centers:
\begin{equation}\label{eq:indexSetsCells}
    \begin{aligned}
    \mathcal{F}^c &= \{i : \mathcal{V}_i \in \mathcal{V}, \quad  \mbox{and a fine time-step $\Delta t / M$ is used on $\mathcal{V}_i$.}  \},\\
    \mathcal{I}^c_1 &= \{i : \mathcal{V}_i \in \mathcal{V}, \quad  \mbox{and $\mathcal{V}_i$ shares at least one edge with some $\mathcal{V}_j$, with $j\in \mathcal{F}^c$.}  \},\\
        \mathcal{I}^c_2&= \{i : \mathcal{V}_i \in \mathcal{V}, \quad  \mbox{and $\mathcal{V}_i$ shares a least one edge with some $\mathcal{V}_j$, with $j\in \mathcal{I}^c_1$.}  \},\\
        \mathcal{C}^c&=\{i: i \notin (\mathcal{F}^c \cup \mathcal{I}^c_1 \cup \mathcal{I}^c_2)\}.
    \end{aligned}
\end{equation}
Note that the above sets form a disjoint cover of the cell index set of the Voronoi tessellation $\mathcal{V}$, i.e. if $\mathcal{V}_i \in \mathcal{V}$, then $i$ is in exactly one of the sets in Eq. \eqref{eq:indexSetsCells}.
The sets identify the LTS regions, and can be seen sketched in Figure \ref{fig:3}: the index set $\mathcal{F}^c$ refers to the fine LTS region, in blue in the figure, $\mathcal{I}^c_1$ contains the cell indices of interface 1, pink in Figure \ref{fig:3}, $\mathcal{I}^c_2$ contains the indices of the cells on interface 2, yellow in the figure, and finally $\mathcal{C}^c$ is the index set for the cells in the coarse LTS region, portrayed in red.
It is necessary to define a few more index sets for the centers:
\begin{equation}
\begin{aligned}\label{eq:index_sets_cells}
    \mathcal{F}_1^c &= \{i \in \mathcal{F}^c : \mbox{$\mathcal{V}_i$ shares at least on edge with some $\mathcal{V}_j$, with $j\in \mathcal{I}_1^c$.}  \}\\
    \mathcal{F}_2^c &= \{i \in \mathcal{F}^c : \mbox{$\mathcal{V}_i$ shares at least on edge with some $\mathcal{V}_j$, with $j\in \mathcal{F}_1^c$.}  \}\\
    \mathcal{\underline{F}}^c &= \mathcal{F}_1^c \cup \mathcal{F}_2^c.
\end{aligned}
\end{equation}
Note that $\mathcal{\underline{F}}^c \subset \mathcal{F}^c$, and it corresponds to the two layers of cells in the fine region that are closest to the interface 1 region.
This set of cells is depicted in white in Figure \ref{fig:3}.
\begin{remark}
Note that the index set $ \mathcal{\underline{F}}^c$ is not needed for the LTS scheme of order two, but it is necessary to guarantee the third order convergence for the LTS scheme of order three. Namely, when advancing the solution on the interface layers with the first stage, one has to also advance the solution on the cells in  $\mathcal{\underline{F}}^c$ with the \it{coarse time-step}. In this way, when computing the second stage on the interface layers, the right-hand side term can be evaluated using the update solution on the fine and not just the old solution as it would be without advancing the fine solution on $\mathcal{\underline{F}}^c$. Note that because the coarse time-step is used on the fine for this step, it might be convenient to select the LTS regions so that the cells in  $\mathcal{\underline{F}}^c$ have about the same size as the coarse cells. This comment is related to the one previously made about the selection of the interface layer regions with respect to the variable resolution region on the mesh.
The need to consider  $\mathcal{\underline{F}}^c$ was actually not explicitly stated in \cite{hoang2019conservative}. Moreover, if one is following the formulation of LTS3 in \cite{hoang2019conservative}, a similar set should be considered for the coarse region. The reason is that in \cite{hoang2019conservative} it is advised to advance the interface layers with the first and second stage sequentially.
In our formulation here, at the first stage we advance {\it both} the coarse region and the interface layers, hence we do not need an analogous of $\mathcal{\underline{F}}^c$ on the coarse.
\end{remark}

Similar index sets as in \eqref{eq:index_sets_cells} are defined for the edges:
\begin{equation}\label{eq:indexSetsEdges}
    \begin{aligned}
    \mathcal{F}^e &= \{e : \mbox{$e$ is an edge of a cell indexed by $\mathcal{F}^c \setminus \underline{\mathcal{F}}^c$.} \},\\
        \underline{\mathcal{F}}^e &= \{e : \mbox{$e$ is an edge of a cell indexed by $\underline{\mathcal{F}}^c$.} \},\setminus \mathcal{F}^e,\\
    \mathcal{I}^e_1 &= \{e :  \mbox{$e$ is an edge of a cell indexed by $\mathcal{I}^c_1$.}  \} \setminus \underline{\mathcal{F}}^e,\\
\mathcal{I}^e_2 &= \{e :  \mbox{$e$ is an edge of a cell indexed by $\mathcal{I}^c_2$.}  \} \setminus \mathcal{I}^e_1,\\
\mathcal{C}^e &= \{e :  \mbox{$e$ is an edge of a cell indexed by $\mathcal{C}^c$.}  \} \setminus \mathcal{I}^e_2.\\
    \end{aligned}
\end{equation}
The above definition implies that the edges are assigned to the LTS regions moving from the fine to the coarse, as it can be seen from Figure \ref{fig:3}. Namely, edges of cells shared between fine cells and interface 1 cells are assigned to the fine region, edges shared between interface 1 and interface 2 cells are assigned to the interface 1 region, and finally edges shared between interface 2 and coarse region are assigned to the interface 2 region.
%%%%%%%%%%%%
%%%%%%%%%%%%

The idea of the LTS algorithms in \cite{hoang2019conservative} is the following: first the solution at the interface layers is advanced using the coarse time-step, and its values at intermediate fine time levels are predicted using interpolation. This is called the interface prediction step.
Then, the predicted values are used to advance with the fine time-step the solution on the fine regions.
Simultaneously, at least from a theoretical point of view, the solution is advanced on the coarse region with the coarse time-step. Then, the solution at the interface is corrected using the solution advanced on the other regions.
In a practical implementation however, it is not necessary to rigorously follow this recipe and the steps should be rearranged to obtain better efficiency. 
\begin{remark} 
Whereas from a theoretical point of view the advancement of the solution in the fine and in the coarse LTS regions may be done in parallel, in practice the situation is more complex due to the necessity of handling parallel communication among processors. In fact, if ghost cells layers (or halos) are used to handle the parallel communication, it is often not possible to advance the solution in parallel in the fine and the coarse, and appropriate strategies have to be adopted to still achieve proper load balancing and parallel scaling. These strategies will be discussed in detail in the numerical results section.
\end{remark}

We now report the LTS algorithms from \cite{hoang2019conservative} as the authors of this manuscript have implemented them in MPAS. The second order LTS algorithm, from now on referred to as LTS2, is a two-stage method based on SSPRK2, whereas the third order LTS algorithm is a three-stage method based on SSPRK3.
\begin{remark}
A very convenient feature of the LTS algorithms in \cite{hoang2019conservative} is that if the parameter $M \in \mathbb{N}^+$ that defines the fine time-step $\Delta t /M$ is equal to one, the LTS algorithms recover the SSPRK methods on which they are based, i.e. if $M=1$, LTS2 is the same as SSPRK2 and LTS3 is the same as  SSPRK3. This feature of the LTS algorithms facilitates  testing their correct implementation and expected convergence order, as the results obtained with one SSPRK method have to match up to machine precision with those obtained with the associated LTS scheme. 
\end{remark}

\begin{remark}
The presentation of the LTS algorithms given in \cite{hoang2019conservative} is more linear and motivated by mathematical considerations, whereas the one adopted here by the authors is less straightforward but is chosen  to reflect the way the algorithms have been implemented, with special consideration towards minimizing storage requirements and computational cost. For instance, the first stage is performed in this work for both the interface layers and the coarse region, so that one can avoid computing on the two coarse layers closest to interface layer 2 only to discard that computation and do it again when the coarse advancement should be done (which is after the second stage computation on the interface, according to \cite{hoang2019conservative}).   
\end{remark}

 As we did for SSPRK, we are going to present the LTS methods in a compact form as follows. 
\rev{Note that the LTS algorithm below is also summarized in Figure \ref{fig:flow-chart}.
%%%%%%%%%%%%
\begin{figure}[!t]
   \centering
   \includegraphics[scale=0.6]{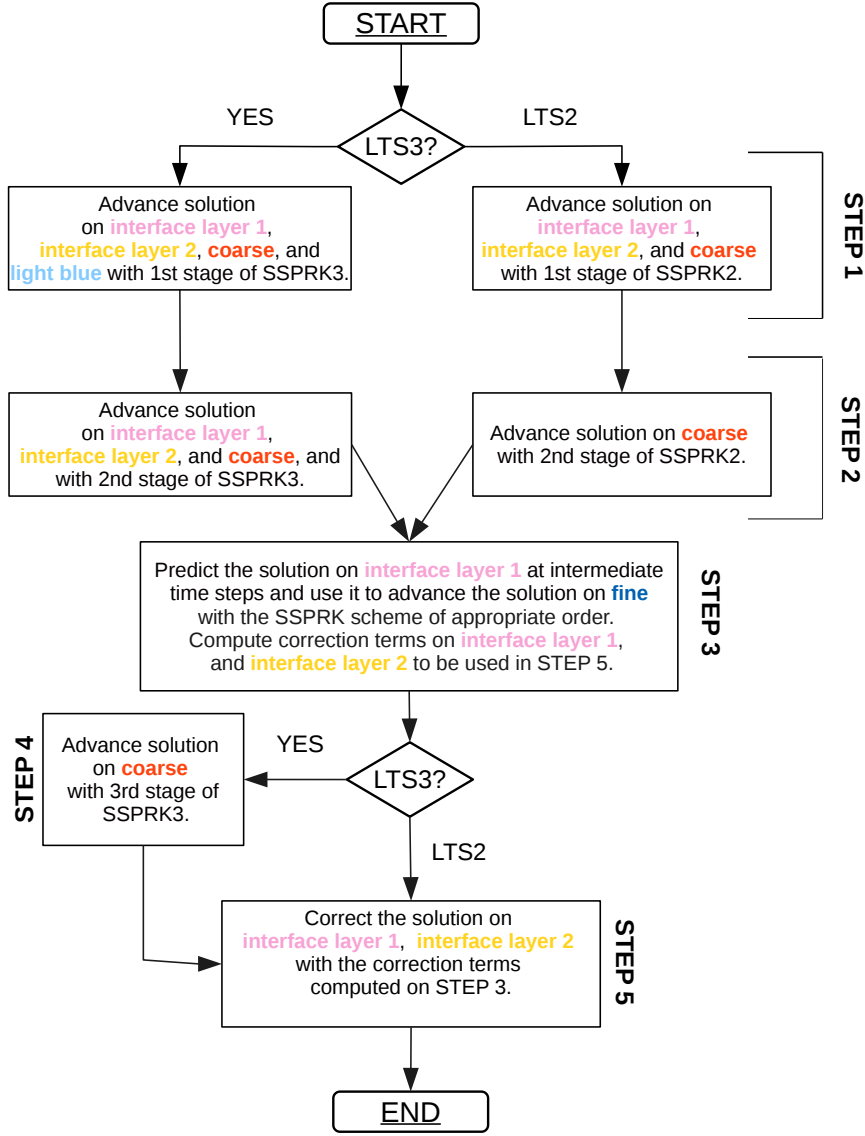}
\caption{Flow chart for Algorithm 2.2. Note that the colors mentioned in the diagram refer to Figure \ref{fig:3}.}
   \label{fig:flow-chart}
\end{figure}
%%%%%%%%%%%%
}

\begin{algorithm}[LTS]Assume $(\mathbf{h}^{old},\mathbf{u}^{old})$, $\Delta t$ (the coarse time-step) and $M \in $ are given. Recall that $M \in \mathbb{N}^+$ is the parameter that defines the fine time-step $\Delta t / M$ that is used to advance the solution on the fine region cells. 
Let us initialize $(\mathbf{h}^{1st},\mathbf{u}^{1st})$ and $(\mathbf{h}^{2nd},\mathbf{u}^{2nd})$ as 
$$(\mathbf{h}^{1st},\mathbf{u}^{1st}) = (\mathbf{h}^{old},\mathbf{u}^{old}) \quad \mbox{and} \quad (\mathbf{h}^{2nd},\mathbf{u}^{2nd}) = (\mathbf{h}^{old},\mathbf{u}^{old}).$$
Then, depending on whether LTS2 or LTS3 is used, $(\mathbf{h}^{new},\mathbf{u}^{new})$ is computed with the following steps:

\underline{Step 1}. Advance the coarse, interface 1 and interface 2 regions with the first stage of the associated SSPRK scheme. Note that if LTS3, then also the two layers of fine cells closest to interface 1 have to be advanced in order to guarantee third order convergence. These layers are indexed by $\mathcal{\underline{F}}^c$ for the cells and $\mathcal{\underline{F}}^e$ for the edges. 
\begin{equation}
\begin{aligned}
\begin{cases}
h_i^{1st} = h_i^{old} + \Delta t \, \mathcal{H}_i(\mathbf{h}^{old},\mathbf{u}^{old})\\
u_e^{1st} = u_e^{old} + \Delta t \, \mathcal{U}_e(\mathbf{h}^{old},\mathbf{u}^{old})
\end{cases}
\end{aligned}
\end{equation}

\begin{equation}
\mbox{for} \quad i \in 
    \begin{cases}
\mathcal{I}_1^c \cup \mathcal{I}_2^c \cup \mathcal{C}^c, \quad {\mbox if LTS2.}\\
\mathcal{\underline{F}}^c \cup \mathcal{I}_1^c \cup \mathcal{I}_2^c \cup \mathcal{C}^c, \quad {\mbox if LTS3.}\\
    \end{cases}, \qquad 
    \mbox{for} \quad e \in 
    \begin{cases}
\mathcal{I}_1^e \cup \mathcal{I}_2^e \cup \mathcal{C}^e, \quad {\mbox if LTS2.}\\
\mathcal{\underline{F}}^e \cup \mathcal{I}_1^e \cup \mathcal{I}_2^e \cup \mathcal{C}^e, \quad {\mbox if LTS3.}\\
    \end{cases}
\end{equation}

\underline{Step 2}. Advance the coarse cells according to the second stage of the associated SSPRK time-stepping. If LTS3, then the solution is advanced also on interface 1 and interface 2.

\begin{equation}
\begin{cases}
h_i^{2nd} = \omega^{old} h_i^{old} + \omega^{1st} h_i^{1st} + \omega^{rhs} \Delta t \, \mathcal{H}_i(\mathbf{h}^{1st},\mathbf{u}^{1st})\\[0.4cm]
u_e^{2nd} = \omega^{old} u_e^{old} + \omega^{1st} u_e^{1st} + \omega^{rhs} \Delta t \, \mathcal{U}_e(\mathbf{h}^{1st},\mathbf{u}^{1st})
\end{cases},
\end{equation}

\begin{equation}
\mbox{for} \quad i \in 
    \begin{cases}
\mathcal{C}^c, \quad {\mbox if LTS2.}\\
 \mathcal{I}_1^c \cup \mathcal{I}_2^c \cup \mathcal{C}^c, \quad {\mbox if LTS3.}\\
    \end{cases}, \qquad 
    \mbox{for} \quad e \in 
    \begin{cases}
 \mathcal{C}^e, \quad {\mbox if LTS2.}\\
\mathcal{I}_1^e \cup \mathcal{I}_2^e \cup \mathcal{C}^e, \quad {\mbox if LTS3.}\\
    \end{cases},
\end{equation}

\begin{equation}
\omega^{old}=\begin{cases}
1/2 \quad \mbox{if LTS2} \\[0.4cm]
3/4 \quad \mbox{if LTS3}
\end{cases},\quad
\omega^{1st}=\omega^{rhs}=\begin{cases}
1/2 \quad \mbox{if LTS2} \\[0.4cm]
1/4 \quad \mbox{if LTS3}
\end{cases}.
\end{equation}
Note that $\omega^{old}, \omega^{1st}$ and $\omega^{rhs}$ correspond to those defined earlier for SSPRK.

\underline{Step 3}. Advance the solution on the fine region. At the same time we build correction terms that are necessary to perform the interface correction. The choice of accumulating the correction terms here is dictated by the need to minimize storage requirements in the implementation. To this end, let us initialize the following quantities:
\begin{equation}
\begin{aligned}
    h_i^{old,0} = h_i^{old}, \quad h_i^{1st,0} = h_i^{1st}, \quad h_i^{2nd,0} = h_i^{2nd} \quad \mbox{for $i \in \mathcal{F}^c \cup \mathcal{I}_1^c \cup \mathcal{I}_2^c.$}\\ 
        u_e^{old,0} = u_e^{old}, \quad u_e^{1st,0} = u_e^{1st}, \quad u_e^{2nd,0} = u_e^{2nd} \quad \mbox{for $e \in \mathcal{F}^e \cup \mathcal{I}_1^e \cup \mathcal{I}_2^e.$}\\
        H_i^{1st,0}=0, \quad H_i^{2nd,0}=0, \quad H_i^{3rd,0}=0, \quad \mbox{for $i \in \mathcal{I}_1^c \cup \mathcal{I}_2^c.$}\\ 
                U_e^{1st,0}=0, \quad U_e^{2nd,0}=0, \quad U_e^{3rd,0}=0, \quad \mbox{for $e \in \mathcal{I}_1^e \cup \mathcal{I}_2^e.$}\\ 
\end{aligned}
\end{equation}
Now, for $k=0,\ldots,M-1$ do:

\underline{Step 3.1.1}. Predict the values of the solution in the interface 1 region at the current fine time interval
\begin{equation}
\begin{cases}
h_i^{old,k} = (1 - \alpha_k - \widetilde{\alpha}_k)h_i^{old} + (\alpha_k - \widetilde{\alpha}_k) h_i^{1st} + 2 \widetilde{\alpha}_k \,\,h_i^{2nd}, \quad \mbox{for $i \in \mathcal{I}_1^c$}\\
u_e^{old,k} = (1 - \alpha_k - \widetilde{\alpha}_k)u_e^{old} + (\alpha_k - \widetilde{\alpha}_k) u_e^{1st} + 2 \widetilde{\alpha}_k \,\,u_e^{2nd}, \quad \mbox{for $e \in \mathcal{I}_1^e$}
\end{cases} ,
\end{equation}

\begin{equation}
\alpha_k = \dfrac{k}{M},\quad
\widetilde{\alpha}_k = \begin{cases}
0 \quad \mbox{if LTS2} \\
\dfrac{k^2}{M^2} \quad \mbox{if LTS3}
\end{cases}.
\end{equation}

\underline{Step 3.1.2}. Compute the correction terms for the interface regions
\begin{equation}
\begin{cases}
H^{1st,k+1}_i = H^{1st,k}_i + \mathcal{H}_i(\mathbf{h}^{old,k},\mathbf{u}^{old,k}), \quad \mbox{ for $i \in \mathcal{I}_1^c \cup \mathcal{I}_2^c$} \\
U^{1st,k+1}_e = U^{1st,k}_e + \mathcal{U}_e(\mathbf{h}^{old,k},\mathbf{u}^{old,k}), \quad \mbox{ for $e \in \mathcal{I}_1^e \cup \mathcal{I}_2^e$}
\end{cases}.
\end{equation}

\underline{Step 3.1.3}. Advance the solution on the fine region

\begin{equation}
\begin{cases}
h_i^{1st,k} = h_i^{old,k} + \frac{\Delta t}{M} \, \mathcal{H}_i(\mathbf{h}^{old,k},\mathbf{u}^{old,k}), \quad \mbox{ for $i \in \mathcal{F}^c$}\\
u_e^{1st,k} = u_e^{old,k} + \frac{\Delta t}{M} \, \mathcal{U}_e(\mathbf{h}^{old,k},\mathbf{u}^{old,k}), \quad \mbox{for $e \in \mathcal{F}^e$}
\end{cases}.
\end{equation}

\underline{Step 3.2.1}. Predict the values of the solution in the interface 1 region at the current fine time interval
\begin{equation}
\begin{cases}
h_i^{1st,k} = (1 - \beta_k - \widetilde{\beta}_k)h_i^{old} + (\beta_k-\widetilde{\beta}_k) h_i^{1st} + 2 \widetilde{\beta}_k \,h_i^{2nd},  \quad \mbox{for $i \in \mathcal{I}_1^c$}\\
u_e^{1st,k} = (1 - \beta_k - \widetilde{\beta}_k)u_e^{old} + (\beta_k-\widetilde{\beta}_k) u_e^{1st} + 2 \widetilde{\beta}_k \,u_e^{2nd},  \quad \mbox{for $e \in \mathcal{I}_1^e$}
\end{cases}.
\end{equation}

\begin{equation}
\beta_k = \dfrac{k+1}{M},\quad
\widetilde{\beta}_k = \begin{cases}
0 \quad \mbox{if LTS2} \\
\dfrac{k(k+2)}{M^2} \quad \mbox{if LTS3}
\end{cases}.
\end{equation}

\underline{Step 3.2.2}. Compute the correction terms for the interface regions
\begin{equation}
\begin{cases}
H^{2nd,k+1}_i = H^{2nd,k}_i + \mathcal{H}_i(\mathbf{h}^{1st,k},\mathbf{u}^{1st,k}), \quad \mbox{ for $i \in \mathcal{I}_1^c \cup \mathcal{I}_2^c$} \\
U^{2nd,k+1}_e = U^{2nd,k}_e + \mathcal{U}_e(\mathbf{h}^{1st,k},\mathbf{u}^{1st,k}), \quad \mbox{ for $e \in \mathcal{I}_1^e \cup \mathcal{I}_2^e$} \\
\end{cases}.
\end{equation}

\underline{Step 3.2.3}. Advance the solution on the fine region

\begin{equation}
\begin{cases}
h_i^{2nd} = \omega^{old} h_i^{old} + \omega^{1st} h_i^{1st} + \omega^{rhs} \frac{\Delta t}{M} \, \mathcal{H}_i(\mathbf{h}^{1st,k},\mathbf{u}^{1st,k}), \quad \mbox{ for $i \in \mathcal{F}^c$}\\
u_e^{2nd} = \omega^{old} u_e^{old} + \omega^{1st} u_e^{1st} + \omega^{rhs} \frac{\Delta t}{M} \, \mathcal{U}_e(\mathbf{h}^{1st,k},\mathbf{u}^{1st,k}), \quad \mbox{for $e \in \mathcal{F}^e$}
\end{cases},
\end{equation}
and $\omega^{old}, \omega^{1st},$ and $\omega^{rhs}$ defined as in Step 2.

If LTS2:

\begin{equation}
    h_i^{new} = h_i^{2nd}, \quad \mbox{for $i \in \mathcal{F}^c$}, \qquad u_e^{new} = u_e^{2nd}, \quad \mbox{for $e \in \mathcal{F}^e$}. 
\end{equation}

Else if LTS3:

\underline{Step 3.3.1}. Predict the values of the solution in the interface 1 region at the current fine time interval
\begin{equation}
\begin{cases}
h_i^{2nd,k} = (1 - \gamma_k - \widetilde{\gamma}_k)h_i^{old} + (\gamma_k-\widetilde{\gamma}_k) h_i^{1st} + 2 \widetilde{\gamma}_k \,h_i^{2nd},  \quad \mbox{for $i \in \mathcal{I}_1^c$}\\
u_e^{2nd,k} = (1 - \gamma_k - \widetilde{\gamma}_k)u_e^{old} + (\gamma_k-\widetilde{\gamma}_k) u_e^{1st} + 2 \widetilde{\gamma}_k \,u_e^{2nd},  \quad \mbox{for $e \in \mathcal{I}_1^e$}
\end{cases}.
\end{equation}

\begin{equation}
\gamma_k = \dfrac{2k+1}{2M},\quad
\widetilde{\gamma}_k = \dfrac{2k^2+2k+1}{2M^2}. 
\end{equation}

\underline{Step 3.3.2}. Compute the correction terms for the interface regions
\begin{equation}
\begin{cases}
H^{3rd,k+1}_i = H^{3rd,k}_i + \mathcal{H}_i(\mathbf{h}^{2nd,k},\mathbf{u}^{2nd,k}), \quad \mbox{ for $i \in \mathcal{I}_1^c \cup \mathcal{I}_2^c$} \\
U^{3rd,k+1}_e = U^{3rd,k}_e + \mathcal{U}_e(\mathbf{h}^{2nd,k},\mathbf{u}^{2nd,k}), \quad \mbox{ for $e \in \mathcal{I}_1^e \cup \mathcal{I}_2^e$} \\
\end{cases}.
\end{equation}

\underline{Step 3.3.3}. Advance the solution on the fine region

\begin{equation}
\begin{cases}
h_i^{new} = \frac{1}{3} h_i^{old} + \frac{2}{3} h_i^{1st} + \frac{2}{3} \Delta t \, \mathcal{H}_i(\mathbf{h}^{2nd},\mathbf{u}^{2nd}) \quad \mbox{ for $i \in \mathcal{F}^c$}\\
u_e^{new} = \frac{1}{3} u_e^{old} + \frac{2}{3} u_e^{1st} + \frac{2}{3} \Delta t \, \mathcal{U}_e(\mathbf{h}^{2nd},\mathbf{u}^{2nd}), \quad \mbox{for $e \in \mathcal{F}^e$}
\end{cases}.
\end{equation}

End if.

Then, set  $h_i^{old,k+1}=h_i^{new}$ for $i \in \mathcal{F}^c$ and $u_e^{old,k+1}=u_e^{new}$.

End do.

\underline{Step 4}. The new solution on the coarse region is obtained.

\begin{equation}
\mbox{ for $i \in \mathcal{C}^c$}: \qquad
h_i^{new} =
\begin{cases} 
h_i^{2nd}, \quad \mbox{if LTS2}\\
\frac{1}{3} h_i^{old} + \frac{2}{3} h_i^{1st} + \frac{2}{3} \Delta t \, \mathcal{H}_i(\mathbf{h}^{2nd},\mathbf{u}^{2nd}), \quad \mbox{if LTS3}
\end{cases},
\end{equation}

\begin{equation}
\mbox{ for $e \in \mathcal{C}^e$}: \qquad
u_e^{new} =
\begin{cases} 
u_e^{2nd}, \quad \mbox{if LTS2}\\
 \frac{1}{3} u_e^{old} + \frac{2}{3} u_e^{1st} + \frac{2}{3} \Delta t \, \mathcal{U}_e(\mathbf{h}^{2nd},\mathbf{u}^{2nd}), \quad \mbox{if LTS3}
\end{cases}.
\end{equation}

\underline{Step 5}. The values at the interface are corrected.

\begin{equation}
\begin{cases}
h_i^{new} = h^{old}_i + \frac{\Delta t}{M}\Big(\theta^{1st} H_i^{1st,M} + \theta^{2nd} H_i^{2nd,M} + \theta^{3rd} H_i^{3rd,M} \Big), \quad \mbox{for $i \in \mathcal{I}_1^c \cup \mathcal{I}_2^c$}\\
u_e^{new} = u^{old}_e + \frac{\Delta t}{M}\Big(\theta^{1st} U_e^{1st,M} + \theta^{2nd} U_e^{2nd,M} + \theta^{3rd} U_e^{3rd,M} \Big), \quad \mbox{for $e \in \mathcal{I}_1^e \cup \mathcal{I}_2^e$}
\end{cases},
\end{equation}

\begin{equation}
\theta^{1st}=\begin{cases}
1/2 \quad \mbox{if LTS2} \\
1/6 \quad \mbox{if LTS3}
\end{cases}, \quad
\theta^{2nd}=\begin{cases}
1/2 \quad \mbox{if LTS2} \\
1/6 \quad \mbox{if LTS3}
\end{cases}, \quad
\theta^{3rd}=\begin{cases}
0 \quad \mbox{if LTS2} \\
2/3 \quad \mbox{if LTS3}
\end{cases}.
\end{equation}

\end{algorithm}

\section{Numerical Investigation}\label{numRes}

The authors of this manuscript have implemented the LTS algorithms proposed in \cite{hoang2019conservative} in the shallow water core of MPAS.
\rev{Using the expertise gained with the present work, we plan to integrate the methods in the ocean core of MPAS as well in the immediate future.}.

\rev{
\subsection{Considerations on performance}
The focus of our analysis is on computational time and parallel performance.
To this end, the three main features that we desire from our implementation of the LTS algorithms are shorter computational time than existing methods, parallel efficiency, and parallel scalability. Here we discuss each of these goals in detail:
\begin{enumerate}
    \item {\underline{Shorter computational time than existing methods}}: we want our implementation to show a computational gain in terms of time when used on meshes with large ratios between the largest and the smallest cell, and a relatively small number of small cells. The gain we intend to observe is relative to explicit global time-stepping methods, in particular the Runge-Kutta time stepping scheme of order four (RK4).  
    To achieve this, it is key to be able to compute the right-hand side terms in \eqref{eq:sweTimeStepping} precisely on a given LTS region, to avoid computing information that is not used at the time of the solution advancement.
    For instance, if we are advancing on the interface cells, we do not want to compute the right-hand side terms on more cells than the interface layers plus two extra layers of cells (needed for the operators in the right-hand side). While this sounds simple in theory, from a parallel computing standpoint it is far from being a trivial task. Luckily, the MPAS framework comes to the rescue and the use of multiple memory blocks per processor makes this goal possible. Specifically, we proceed in the following way:
    \begin{enumerate}
        \item We construct a {\it{graph.info}} file using the graph partitioning library METIS \cite{karypis1997metis}, dividing the domain according to a load balancing strategy only dictated by the cells sizes (so it is irrespective of the LTS regions).
        This file is just a column of numbers where the row index is the cell index and the number on a given row is the MPI processor the cell is assigned to.
        Note that if we were to stop here, we would select the standard approach to achieve fast computational time through parallelism, i.e. the one that is based only on cell sizes. This is perfectly fine for most algorithms but it is definitely not enough for the LTS schemes that we address. As a matter of fact, a partitioning that takes into account only the cell sizes for load balancing will likely produce MPI subdomains that own only coarse time-step cells or only fine time-step cells, resulting in a situation where a given processor is often not carrying out any computation and so unnecessary wait times arise. Other even worse scenarios are possible, such as for instance one where an MPI processor owns mostly cells in the coarse LTS region and them some cells in the interface layers and/or in the fine region. In such a case any time the time-stepping is carried out on these minority regions, the computation of the right-hand side terms would still be carried out on {\it{all}} the cells owned by the processor, even on those where the solution is not advanced, resulting in a largely unnecessary computational effort that dramatically hinders the computational time performance.
        In general, the goal of computing the right-hand side terms on specific LTS regions is very complex, as one would have to handle the interaction of the MPI halo layers with an independent layer of halo cells around the LTS region necessary for the proper computation of the right-hand side term. Obviously these two sets of halos will not match in general, hence the difficulty of the task. An efficient way to make these halos match and therefore handle the issues is discussed next.
        \item Based on the LTS regions assigned to the cells of the mesh, the {\it{graph.info}} is modified by adding to each row a triplet whose entries are weights that identify the different LTS regions. For instance, (1,0,0) identifies the coarse LTS region, (0,1,0) the fine and (0,0,1) the interface layers.
         We do it by using the multi-constraint formulation from METIS (see \cite{karypis1998multilevel} and also the METIS manual).
        With the {\it{graph.info}} file modified in this way, each MPI partition now has some cells from the fine LTS region, some cells from the coarse LTS region, and some from the interface layers, in a balanced way. This means that all processors have approximately the same number of cells in total. Note that this does not imply that the number of fine, interface and coarse cells will be balanced within a given MPI rank, i.e. unless extra cells are added to the interface layers as it will be discussed next, the interface layer cells will be considerably less than the coarse and fine cells within a given MPI partition.
        Note that after this step, we have produced an MPI partitioning where all MPI processors will be working during the time-stepping procedure, because they all possess some cells from all LTS regions. However, while this is an important step towards proper scalability, it is absolutely not sufficient to make the LTS algorithms performant and faster than existing explicit global time-stepping methods because at this point the computation of the right-hand side terms is still carried out on all cells of the MPI partition, hence on all LTS regions owned by the processor, regardless of whether the computation will then be used to advance the solution.
        \item Using the {\it{graph.info}} file defined above, we produce a {\it{graph.info.part.X}} file, where $X=3N$ is the total number of blocks we are going to use, and $N$ is the number of MPI cores. Note that we are using 3 blocks per core because we have 3 different LTS regions, i.e. fine, coarse and interface.
        Then, given that for every processor index $k=0,\ldots,N-1$, the indices of its blocks will be $3k+b$, where $b=0,1,2$, we modify the {\it{graph.info.part.X}} file so that all the first blocks of each processor only own fine cells, all the second blocks own coarse cells and all the third blocks own interface cells (this correspondence is arbitrary and was chosen for ease of implementation).
        Then we run the simulation using the {\it{graph.info.part.X}} file but requesting $N=X/3$ MPI cores. In this way, the simulation runs with $N$ MPI cores, each of which has three blocks, and each block is made of exactly one LTS region. The beauty of this approach is that, with appropriate flags, one can compute the right-hand side terms only on a single specific LTS region plus two halo layers (which are needed by the definition of the operators) which now automatically coincide with the cells own by the blocks, plus their MPI halos. 
        No extra work is necessary to handle the communication and interplay of these two types of halo layers because now they match exactly. 
        Then, on the time stepping scheme, we advance the solution only on the cells and edges {\it{owned}} by a given block (so this time excluding the halo) and afterwards perform an update of the MPI halo layer on the solution.
        We remark that while the blocks have their own MPI halos that overlap with the overall MPI halo of the MPI partition that owns them, they do not take part in the parallelism, as the MPI processor loops over the blocks in a serial manner. This is why it is important that each MPI partition has three blocks, each owning a different LTS regions, so that work can be performed at all times by a given processor without wait times.
        Note that one could use the {\it{graph.info.part.X}} file to run with $X$ MPI cores as well, however in this case each MPI processor will be composed of a single block made up of a single LTS region, hence wait times will not be balanced so it is not the most efficient setup. 
            \end{enumerate}
            \item {\underline{Parallel efficiency}}: the load has to be properly balanced among MPI processors. Note that while with other time-stepping schemes this can be easily handled with METIS only requesting that the load be balanced according the the cell sizes, for the case of LTS this is not straightforward, as the load of a given processor depends on the type of LTS regions that are assigned to it. For instance, if a processor only owns coarse cells, then during the sub-stepping procedure it will have to wait without performing any work. This is not optimal from a load balancing point of view and it has to be addressed, and we do it as explained in point (b) of the previous bullet point.
            In this way, every MPI processor has some cells in the coarse LTS region, some cells in the fine LTS region, and some in the interface, in a load balanced way, so that all processors are performing work at all times during the simulation, and in total they own approximately the same number of cells. Note that parallel efficiency is important for the first goal, which is that of minimizing the computational time.
            \item {\underline{Parallel scalability}}: we want our compute time to scale linearly against the number of MPI processors. For LTS schemes, this is a tricky goal because, as mentioned above, the load depends on the specific LTS region a given processor handles. 
            Moreover, one cannot advance the solution on the fine and on the coarse simultaneously, hence load has to be distributed according to different criteria, as we have already described above.
            Furthermore, the interface layers will in general be composed of fewer cells than the fine and the coarse regions, hence to achieve scalability with our three-block per MPI rank approach, one has to increase the number of cells in the interface layers, as it will be shown next.
\end{enumerate}
}

\rev{\begin{remark} A generalization to the case of more than two LTS regions (for instance a very high resolution zone, a high resolution and a low resolution zone) can be handled theoretically provided that for any two neighboring regions, the time-step size on one be an integer multiple of the time-step on the other, see Remark 4.3 in \cite{hoang2019conservative}.
Our approach extends to such a case in a natural way by choosing appropriate weights with METIS so that each MPI rank owns cells from all LTS regions in a balanced way (as in the case treated here). Moreover, an extra number of blocks corresponding to the additional number of interfaces and time-step regions should be used for the partitioning. For instance, in case of a very high resolution zone, a high resolution and a low resolution zone, one would run with 3 blocks as in the original approach, plus 2 extra blocks, one for the extra interface between very high resolution and high resolution, and one for the very high resolution region.
\end{remark}
}

The numerical results are obtained considering a standard test cases for the shallow water equations on a sphere, namely test case 5 from  \cite{williamson1992standard}, i.e. a zonal flow over an isolated mountain.
The mountain is placed at a point with longitude $\lambda_c= 3\pi/2$ and latitude $\theta_c = \pi/6$, and has a height $h_s$ that depends on longitude-latitude $(\lambda,\theta)$ location of a point on the sphere:  
\begin{align}\label{eq:mnt}
    h_s(\lambda,\theta) = 2000\Big(1-\dfrac{r(\lambda,\theta)}{R}\Big),
\end{align}
where $R=\pi/9$, and $r^2(\lambda,\theta)=\min\{R^2, (\lambda-\lambda_c)^2+(\theta-\theta_c)^2\}$.
The high resolution region throughout will always be centered at point where the top of the mountain is located.

\subsection{Convergence rate test}
An analytic solution is not known for the test case considered, so to compute errors we will use a reference solution obtained with RK4 with a very small time-step, i.e. $\Delta t_{ref} = 0.1$ seconds.

\begin{figure}[!t]
   \centering
   \includegraphics[scale=0.4]{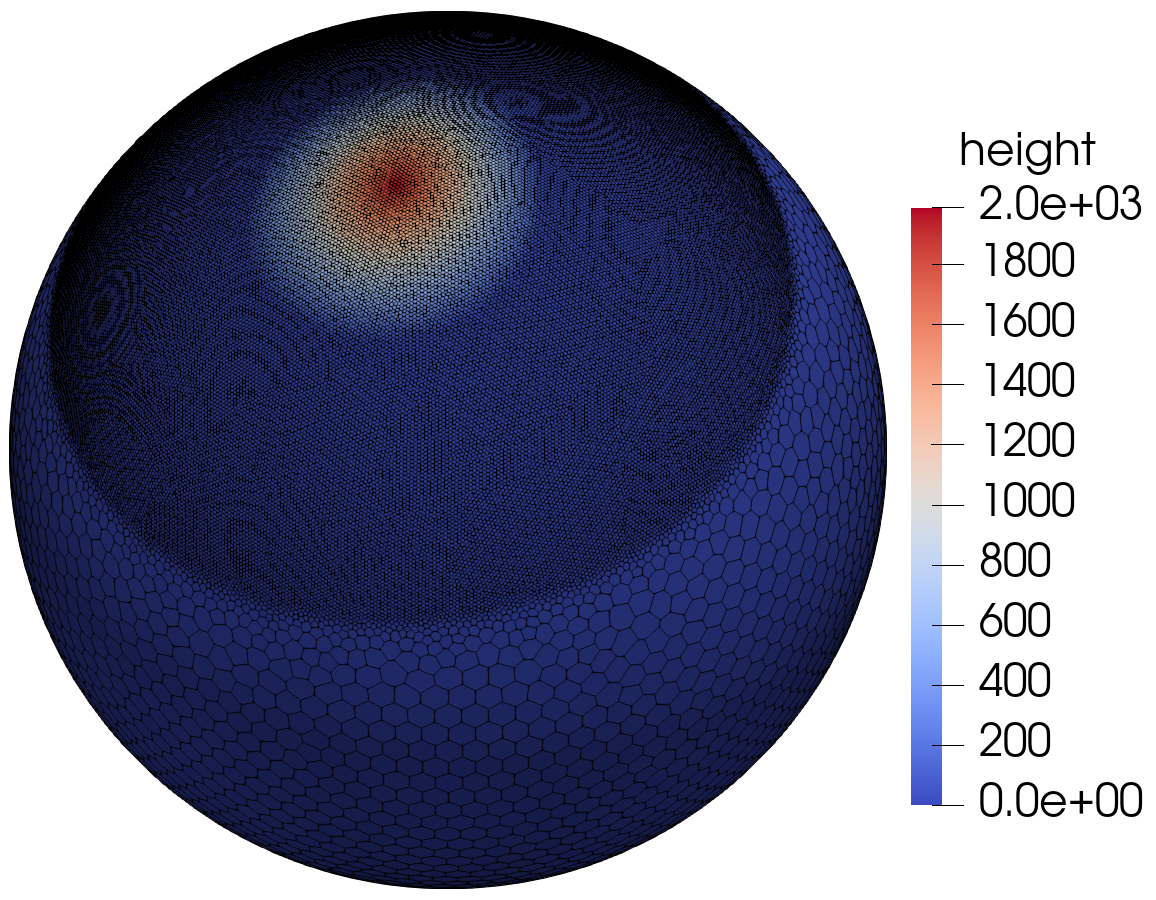}
    \includegraphics[scale=0.4]{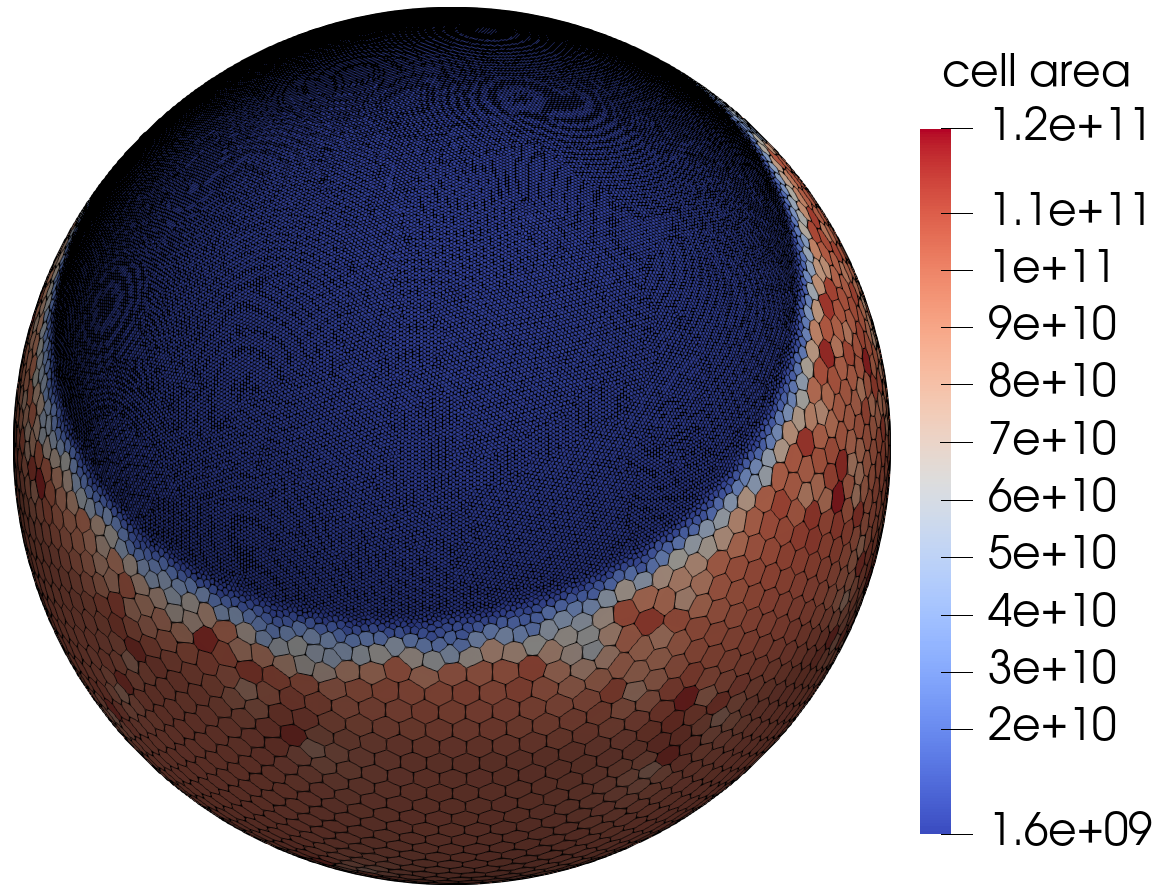}
\caption{Left: profile of the mountain described in Eq. \eqref{eq:mnt}. Right: plot of the cell area of the mesh cells. }
   \label{fig:4}
\end{figure}

\begin{figure}[!t]
   \centering
   \includegraphics[scale=0.4]{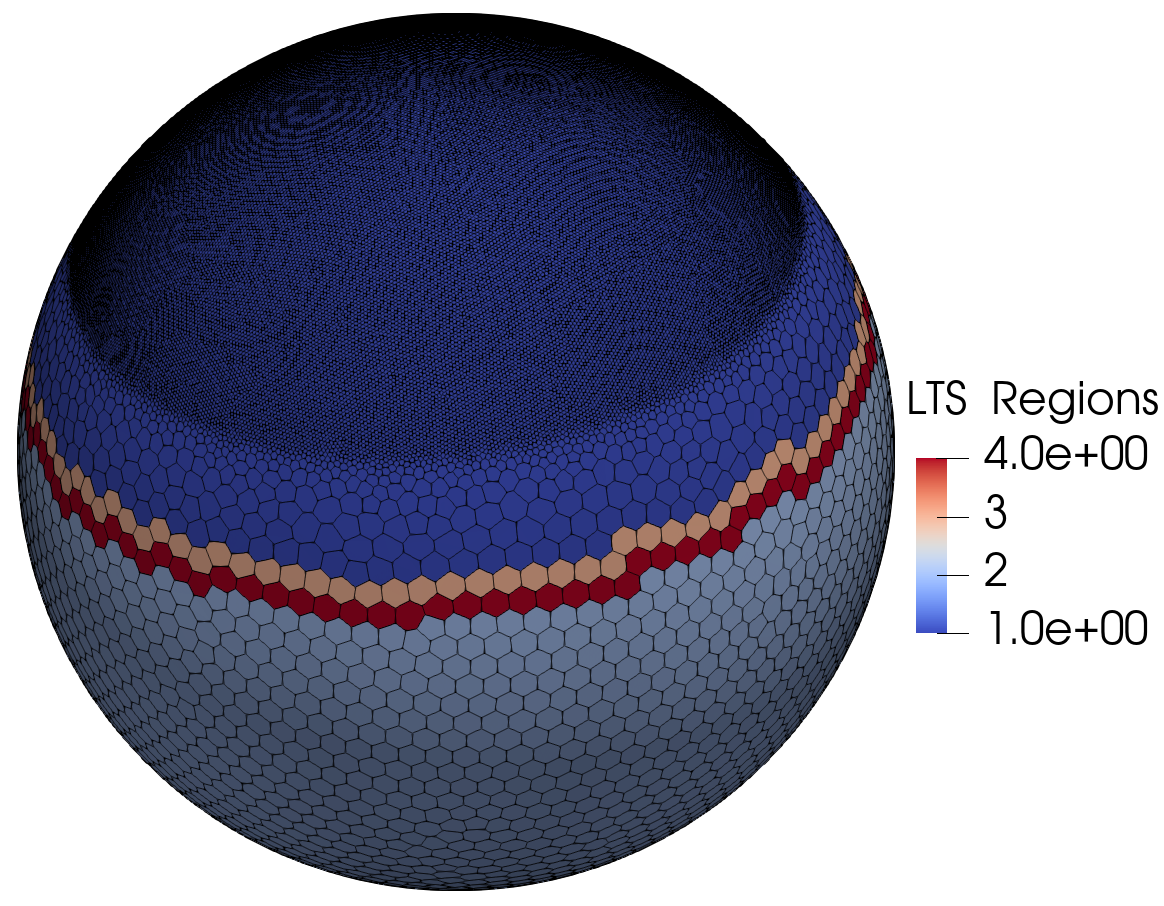}
      \includegraphics[scale=0.4]{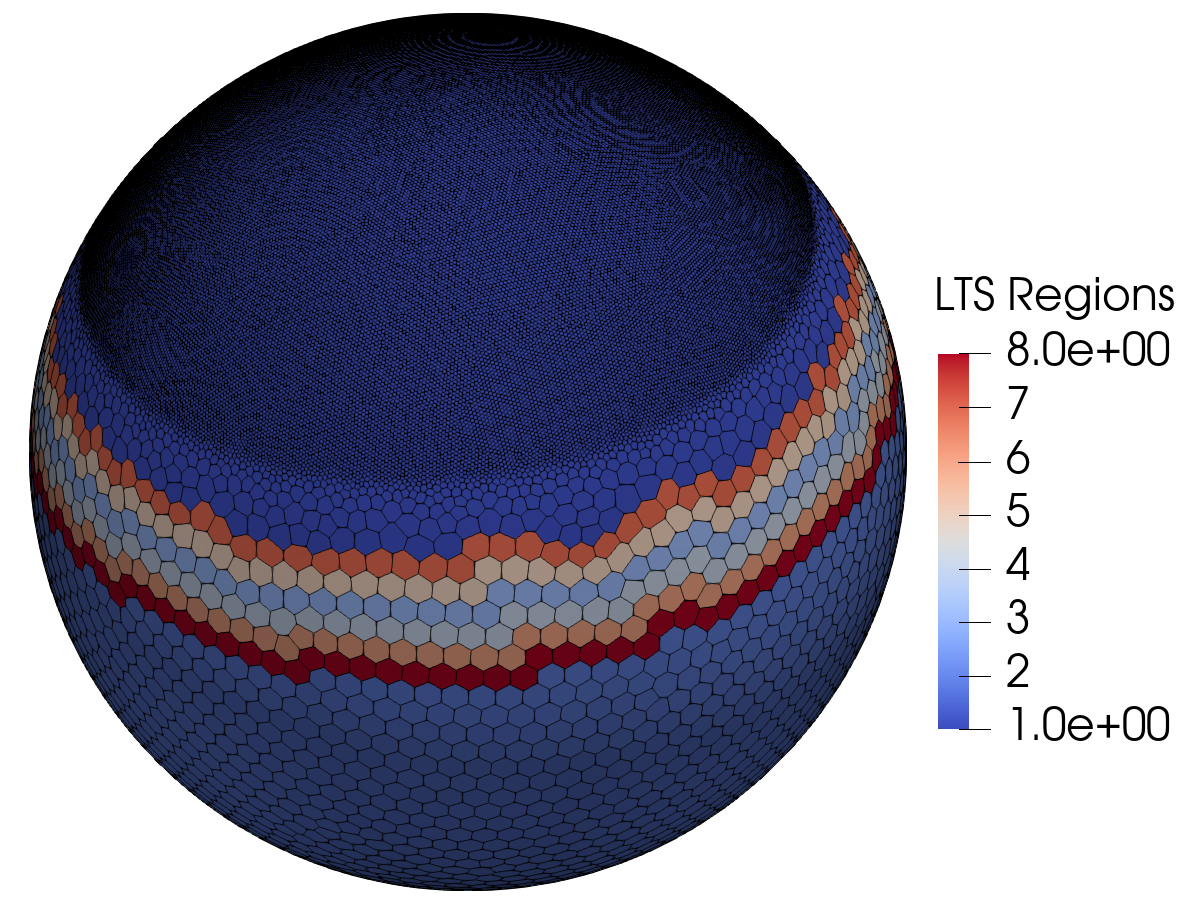}
\caption{Left: LTS regions used for LTS2 to produce the results in Figure \ref{fig:6}. Right: LTS regions used for LTS3 to produce the results in Figure \ref{fig:7}.
LTS Regions legend: fine LTS region = 1,5,7. Coarse LTS region= 2,6,8. Interface layer 1= 3. Interface layer 2= 4.}
   \label{fig:5}
\end{figure}

The height profile of the mountain and the cell areas for the mesh considered in this test are shown in Figure \ref{fig:4}.
Note that the high resolution region completely covers the mountain and also extends around it.
The LTS regions are shown in Figure \ref{fig:5} for LTS2 on the left and LTS3 on the right.
Referring to the legend in Figure \ref{fig:5}, the LTS regions with numbers 1,5 or 7 compose the fine region; numbers 2, 6, or 8 refer to the coarse LTS region, 3 is the interface layer 1 and 4 is the interface layer 2.
Recall that the coarse time-step is used to advance the solution on the interface layers, this is why we chose to locate the interface layers on cells that have the same size as those in the coarse LTS region.
Within LTS3, the cells with an LTS region value of 5 and 7  are those indexed by $\mathcal{\underline{F}}^c$, whereas those with LTS region value of 6 and 8 are treated together with the other coarse cells with LTS region value of 2. Recall that this is possible because in our implementation of the LTS schemes we advance with the first stage the interface layers and the coarse region together. If in LTS3 one was to advance  interface layer 2 with the first {\it and} the second stage and {\it then} advance the coarse, it would be necessary to perform the first stage also on those cells that have an LTS region value of 6 and 8, and then do it only on those with a value of 2 when it is time to advance the coarse. In our approach, for LTS3 we perform the first stage on the interface layers, all the coarse region, and the cells with LTS region 5 and 7. 

%%%%%%%%%%%%
%%%%%%%%%%%%
%%%%%%%%%%%%
%%%%%%%%%%%%

As remarked previously in the manuscript, a convenient feature of the LTS algorithms here considered is that if $M=1$, then the associated SSPRK algorithms are recovered. This is a very helpful tool when it comes to the implementation of the methods, and the first thing one should check is that machine precision matching is obtained between an LTS implementation with $M=1$ and an independent SSPRK implementation of the same order.

\begin{figure}[!t]
   \centering
      \includegraphics[scale=0.51]{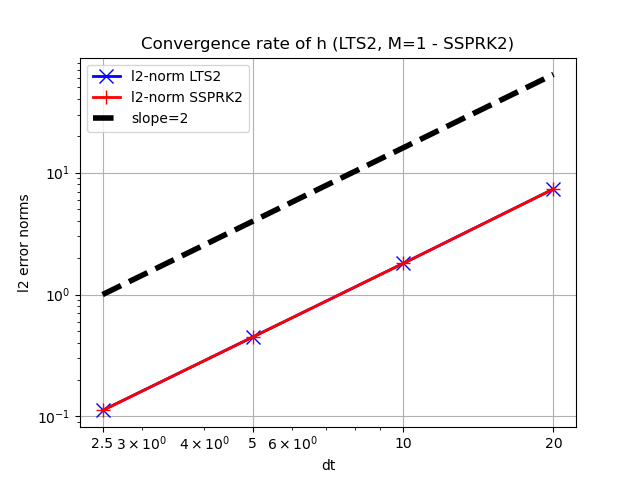}
      \includegraphics[scale=0.51]{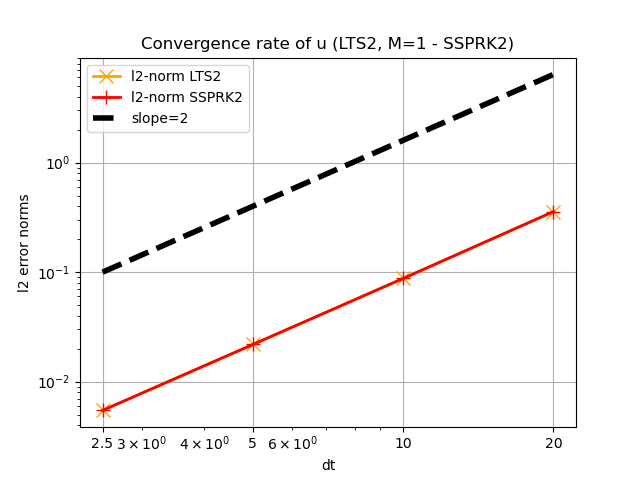}\\
   \includegraphics[scale=0.51]{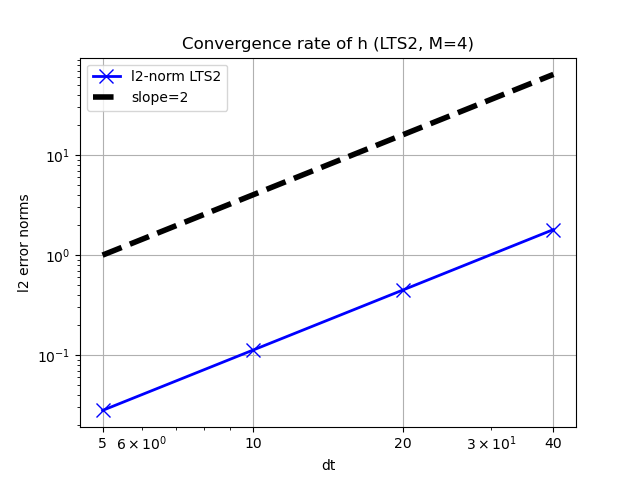}
      \includegraphics[scale=0.51]{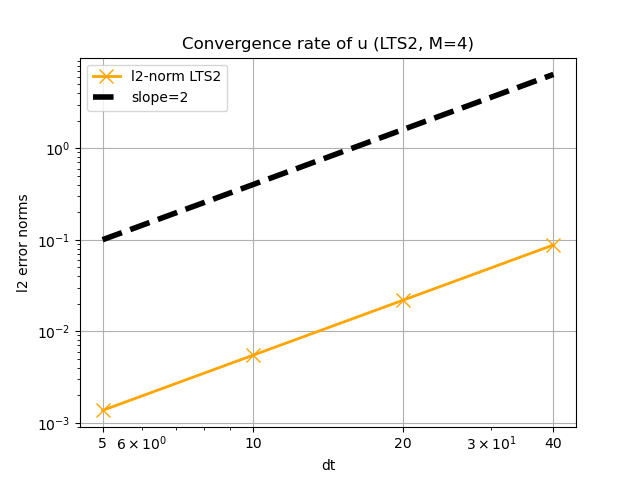}
\caption{Convergence rates for the LTS2 algorithms obtained with respect to a reference RK4 solution with $\Delta t_{ref}=0.1$ seconds. The duration of the simulation is 1 hour. The $\Delta t$ on the x-axis are the coarse time-steps.}
   \label{fig:6}
\end{figure}

\begin{figure}[!t]
   \centering
      \includegraphics[scale=0.51]{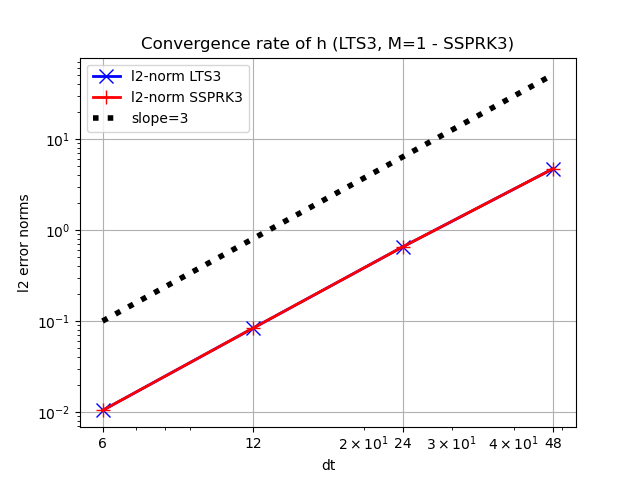}
      \includegraphics[scale=0.51]{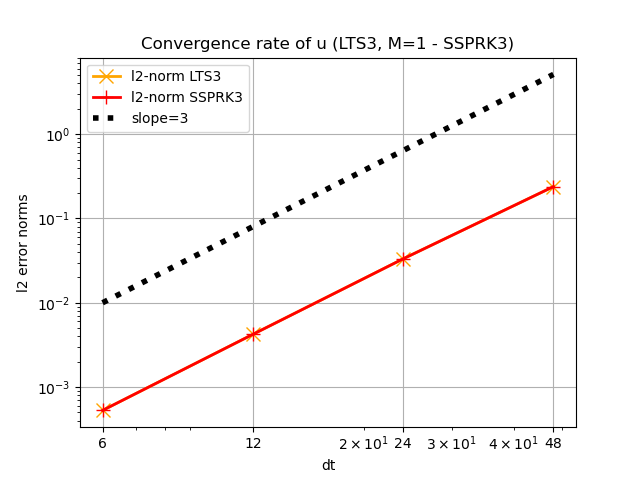}\\
   \includegraphics[scale=0.51]{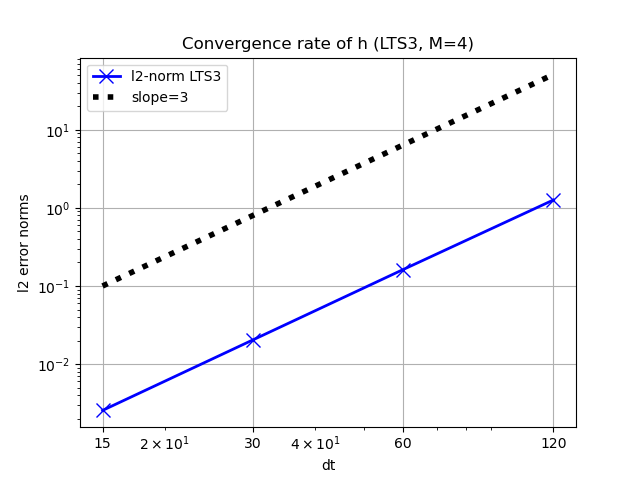}
      \includegraphics[scale=0.51]{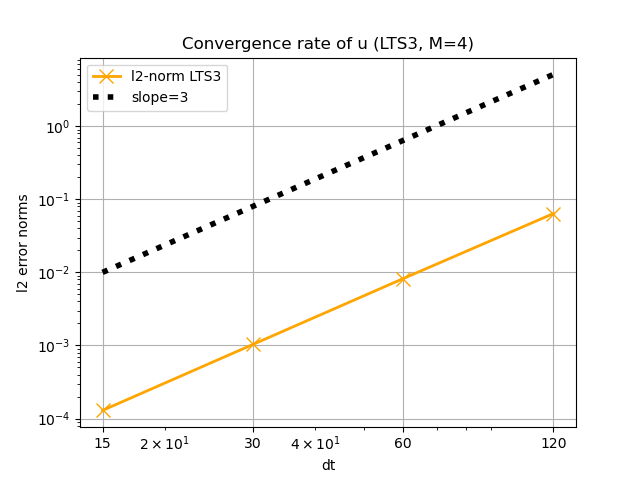}
\caption{Convergence rates for the LTS3 algorithms obtained with respect to a reference RK4 solution with $\Delta t_{ref}=0.1$ seconds.
The duration of the simulation is 1 hour. The $\Delta t$ on the x-axis are the coarse time-steps.}
   \label{fig:7}
\end{figure}

To show that the implemented LTS schemes possess the expected convergence orders we compute 
$$\| sol_{LTS}(\Delta t) - sol_{RK4}(\Delta t_{ref})\|_{l_2},$$ 
for a run of 1 hour, where $sol$ is either the thickness $h$ or the normal velocity $u$. The LTS schemes are used with $M=1$ and $M=4$. For LTS2, with $M=1$ the coarse time-steps (same as the fine) are $\Delta t = 20, 10, 5, 2.5$, with $M=4$ they are $\Delta t = 40, 20, 10, 5$, and hence the fine-time steps will be $10, 5, 2.5, 1.25$.
For LTS3, with $M=1$ the coarse time-steps (same as the fine) are $\Delta t = 48, 24, 12, 6$, whereas with $M=4$ they are $\Delta t = 120, 60, 30, 15$ and hence the fine-time steps will be $30, 15, 7.5, 3.75$.
Results are shown in Figure \ref{fig:6} for LTS2 and Figure \ref{fig:7} for LTS3. Results for independent implementations of SSPRK2 and SSPRK3 are also included. The expected convergence rates are obtained. Note that the plots also show that the LTS algorithms produce the same errors as the associated SSPRK algorithms if $M=1$, up to machine precision.

%%%%%%%%%%%%
%%%%%%%%%%%%
%%%%%%%%%%%%
%%%%%%%%%%%%

\subsection{Gain of local time-stepping and parallel scalability}

\begin{figure}[!t]
   \centering
   \includegraphics[scale=0.3]{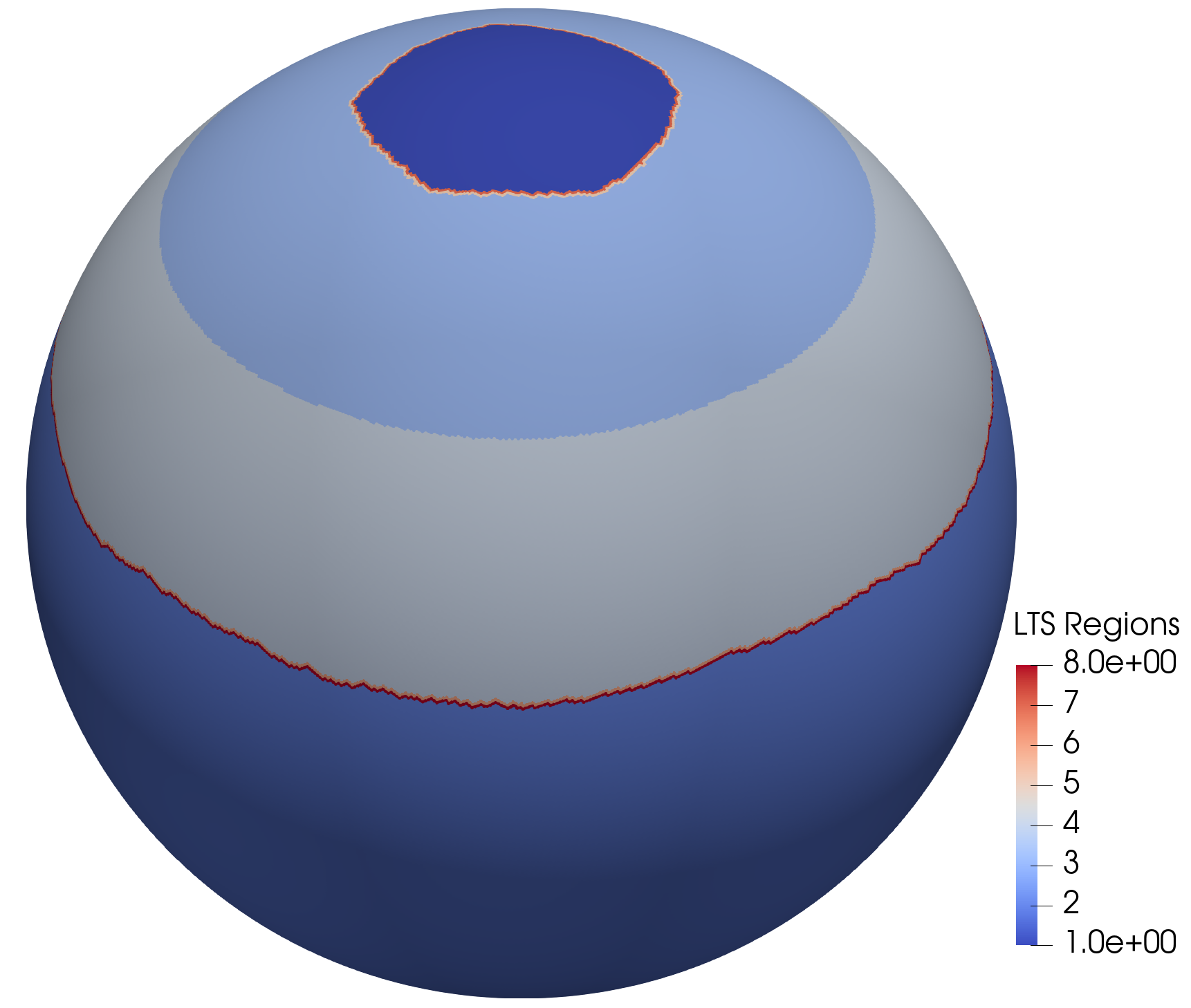}
\caption{LTS regions used for the LTS scheme of order 3 (LTS3) to produce the results in Table \ref{table:2}.  The two interface layers (interface layer 1, which has LTS Region = 3, and interface layer 2 with LTS Region =4) are composed of 78 layers of cells each. This is necessary to achieve parallel scalability.}
   \label{fig:8}
\end{figure}

Next, we report a similar test as what was done in \cite{hoang2019conservative} to show the extent of the computational gain in terms of time, for a given fine time-step, that can be obtained with local time stepping.
This test compares the CPU time of runs with $M=1$ and $M=4$ for a given LTS algorithm and given time-step. Here we use a higher resolution mesh than the one used before, with a total of $621007$ cells, see Figure \ref{fig:8} where the LTS regions are shown for LTS3. Note that for this test the interface layers are not composed of just one layer of cells, but they are actually made up of 78 layers of cells each, see the light blue (interface layer 1) and grey (interface layer 2) regions in the figure.
The reasons for this is that, as we mentioned, for a given MPI rank, we are using 3 blocks, each of which has only one LTS region. This means that we need enough cells in the interface so that, when they are split among MPI ranks, there are enough to compensate the communication time and we can see proper parallel scalability. 

\begin{figure}[!t]
   \centering
   \includegraphics[scale=0.52]{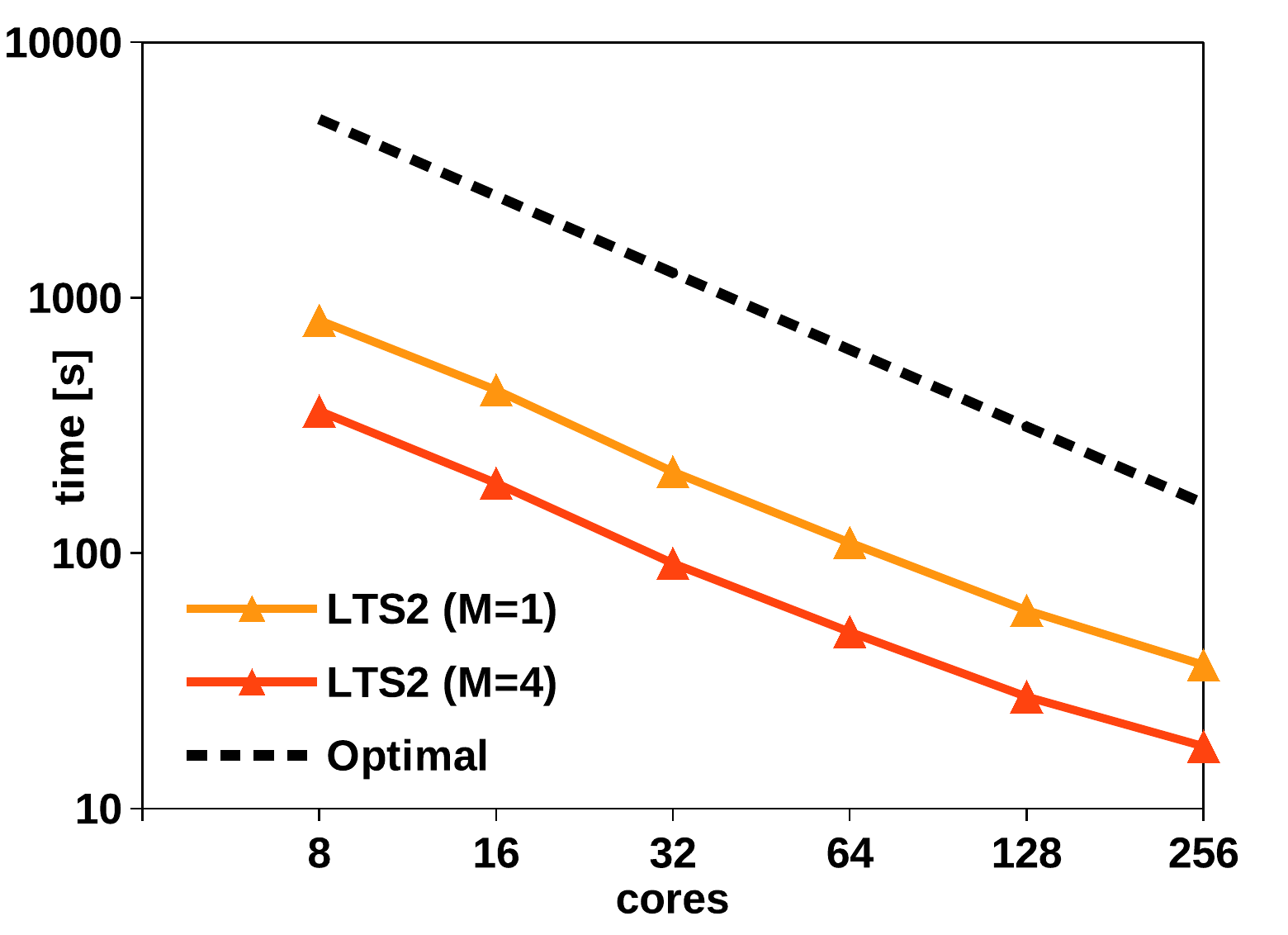}
   \includegraphics[scale=0.52]{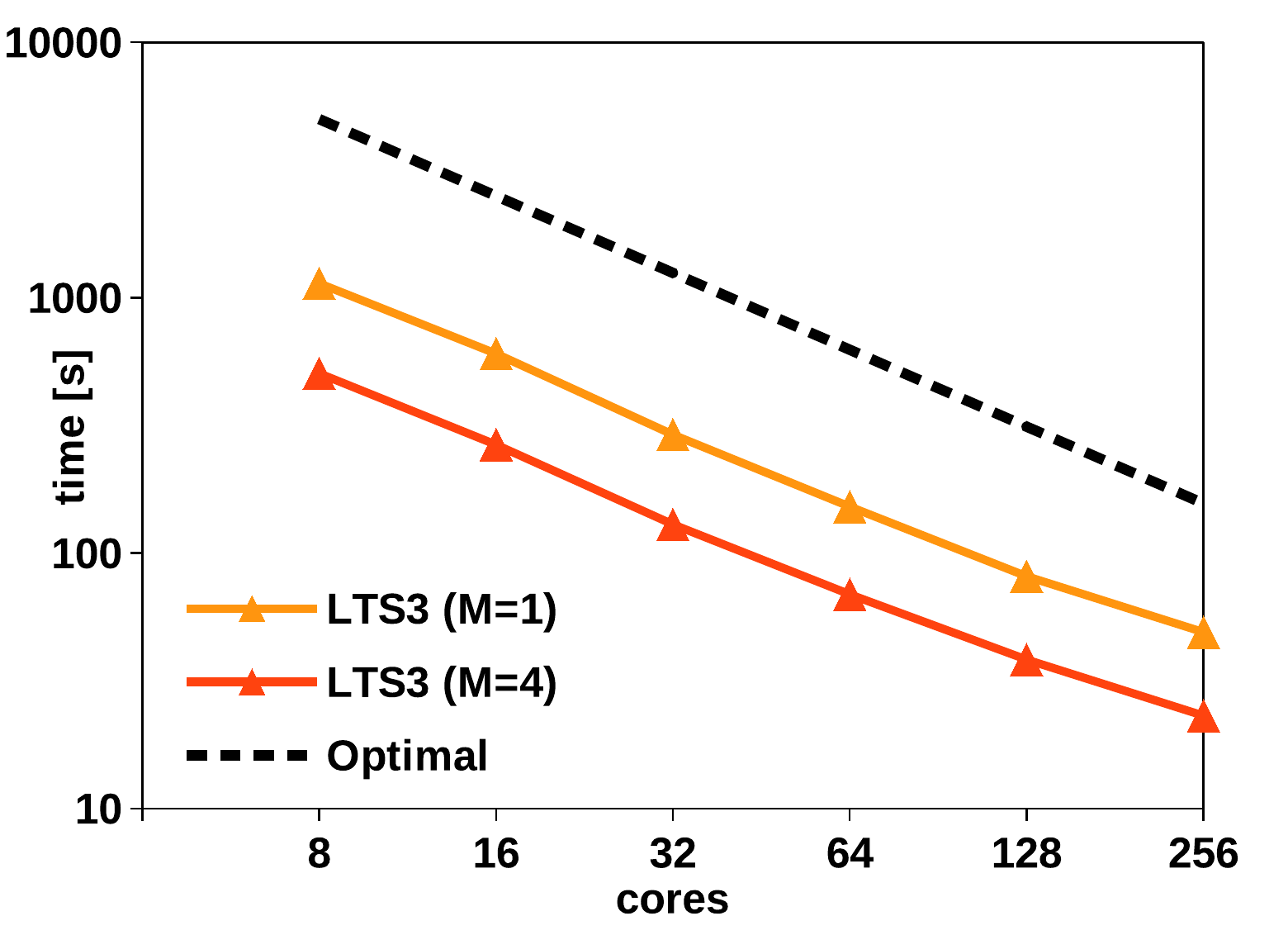}
\caption{Parallel scalability plots. Left: LTS2, the CPU times are in Table \ref{table:1}. Right: LTS3, the CPU times are in Table \ref{table:2}.}
   \label{fig:9}
\end{figure}

The simulations are run in parallel with $N=8,16,32,64,128$ and $256$ processors. Recall that, when running the LTS schemes, because we are using 3 blocks per MPI processor, we need to create with METIS the following graph.info files: graphinfo.file.part.$X$, with $X=3N=24,48,96,192,384,768$. 
Results for LTS2 are reported in Table \ref{table:1}, whereas results for LTS3 are in Table \ref{table:2}. Just as a reference, we also include results for RK4 as wells as the SSPRK schemes associated to the LTS schemes. \rev{The simulations have been run for 2 hours with a fine time-step of $\Delta t = 3$ seconds for all methods. When $M=4$, the coarse time-step is 12 seconds.} We consider only the CPU time required for the integration procedure, and neglect all other times. The initialization time for the LTS regions is negligible compared to the integration time. 
The optimal ratio for the LTS gain is defined as ($4 \, \times$ total number of cells) / (coarse time-step cells + $4 \, \times$ fine time-step cells).
This is the ratio between the number of cells on which the fine time step is used with no local time stepping (which is all cells) divided by where it is used with local time stepping. 
The coarse time-step cells are $474467$ (350816 coarse and 123651 in total for the two interface layers) and the fine time-step cells are $146540$, so the optimal ratio for the mesh considered is approximately $2.342$. For both methods, the speed-up achieved by the local time-stepping procedure is very close to the theoretical one, see Table \ref{table:3}.
Also, the LTS schemes with $M=4$ are the fastest methods, for the time-step considered.
Note that for SSPRK and RK4 the MPI partitions are obtained with METIS only balancing the number of cells per processor, irrespective of the LTS regions. We believe this is the fair approach to adopt, given that RK4 and SSPRK should not see the LTS regions, and so the load balancing for them should only be dictated by the cells size of the mesh.
In Figure \ref{fig:9} we report a plot to visualize the parallel scalability of the LTS algorithms compared to an optimal linear rate. The scalability is close to linear and only gets slightly worse as the number of cores becomes large, recall that our LTS implementation needs 3 times as many blocks as the number of cores.

%%%%%%%%%%%%
%%%%%%%%%%%%

\begin {table}[!h!]
\caption{\rev{Computational times in seconds for LTS2 with $M=1$ and LTS2 with $M=4$. The fine time-step is $\Delta t = 3$ seconds and the coarse is 12 seconds. The run has a duration of 2 hours.}}
\setlength\tabcolsep{3.25pt} % default value: 6pt
\begin{center}
\begin{tabular}{|c|c|c|c|c|} 
\cline{1-5}    
 & \multicolumn{4}{c|}{CPU time [s]}  \\ \hline  
   MPI ranks & RK4 & SSPRK2 & LTS2 ($M=1$) & LTS2 ($M=4$) \\ \hline 
        8 &  1294.602  &   647.973  &   815.194   &  360.165  \\ \hline
        16 & 716.714   &   355.957  & 434.915    &  187.761  \\ \hline
     32 &  347.345  &  170.706   &  208.008   &  91.043   \\ \hline
64 &  160.543  &  82.920   &  109.712   &  49.057    \\ \hline
128 &  78.452   & 42.690    &  59.649   &  27.367   \\ \hline
256 &  38.930  &  23.703   &  36.423    &  17.491   \\ \hline
\end{tabular}
\end{center}
\label{table:1}
\end{table}

\begin {table}[!h!]
\caption{\rev{Computational times in seconds for LTS3 with $M=1$ and LTS3 with $M=4$. The fine time-step is $\Delta t = 3$ seconds and the coarse is 12 seconds. The run has a duration of 2 hours.}}
\setlength\tabcolsep{3.25pt} % default value: 6pt
\begin{center}
\begin{tabular}{|c|c|c|c|c|} 
\cline{1-5}    
 & \multicolumn{4}{c|}{CPU time [s]}  \\ \hline  
   MPI ranks & RK4 & SSPRK3 & LTS3 ($M=1$) & LTS3 ($M=4$) \\ \hline
   8 &  1294.602   &   924.768   &  1137.802    &  505.475  \\\hline
   16 &  716.714    &  508.954    &  603.503    &  265.617  \\\hline
     32 &  347.345  &  241.859   &  290.958   & 129.393  \\\hline
64 & 160.543   &   115.299  &   151.906  &  68.862  \\ \hline
128 & 78.452   &  58.972   &  81.073   &  38.246  \\ \hline
256 & 38.930   &  33.046  &  48.980   &  23.066  \\ \hline
\end{tabular}
\end{center}
\label{table:2}
\end{table}

\begin {table}[!h!]
\caption{\rev{Computational gain for the LTS schemes provided by the local time stepping. The first two columns are the ratio of the CPU times reported in Table \ref{table:1} and \ref{table:2}, respectively.}}
\setlength\tabcolsep{3.25pt} % default value: 6pt
\begin{center}
\begin{tabular}{|c|c|c|c|} 
\cline{1-4}  
    MPI ranks & LTS2 ($M=1$) / LTS2 ($M=4$) &   LTS3 ($M=1$)   / LTS3 ($M=4$)   & Optimal Rate\\ \hline 
     8 &   2.263   &  2.251 & 2.342 \\ \hline
     16 &  2.316   &  2.272  & - \\ \hline
     32 &  2.285   &  2.249  & - \\ \hline
     64 &  2.236   &  2.206 &  - \\ \hline
     128 & 2.179   &  2.120  & -\\ \hline
     256 & 2.082   &  2.123 &  - \\ \hline
\end{tabular}
\end{center}
\label{table:3}
\end{table}

%%%%%%%%%%%%
%%%%%%%%%%%%

\subsection{Maximum time-step and CPU time test}
\begin{figure}[!htb]
   \centering
   \includegraphics[scale=0.35]{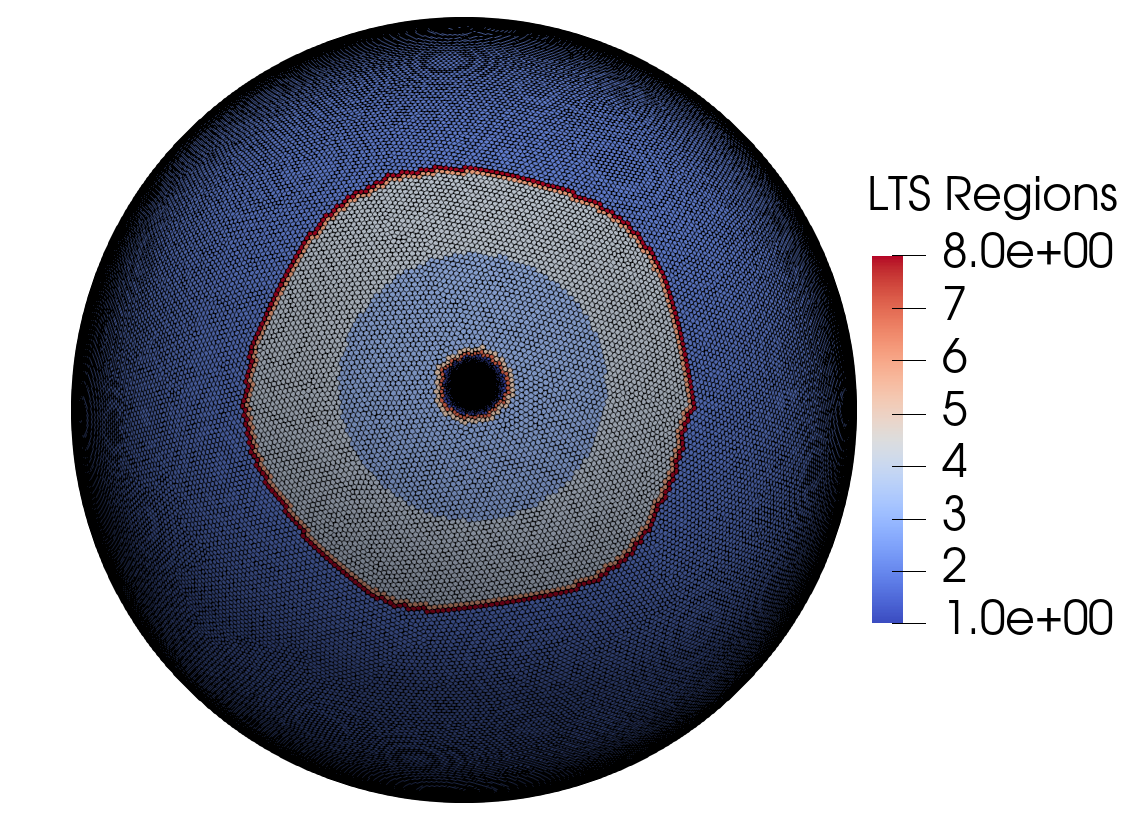}\\
   \includegraphics[scale=0.35]{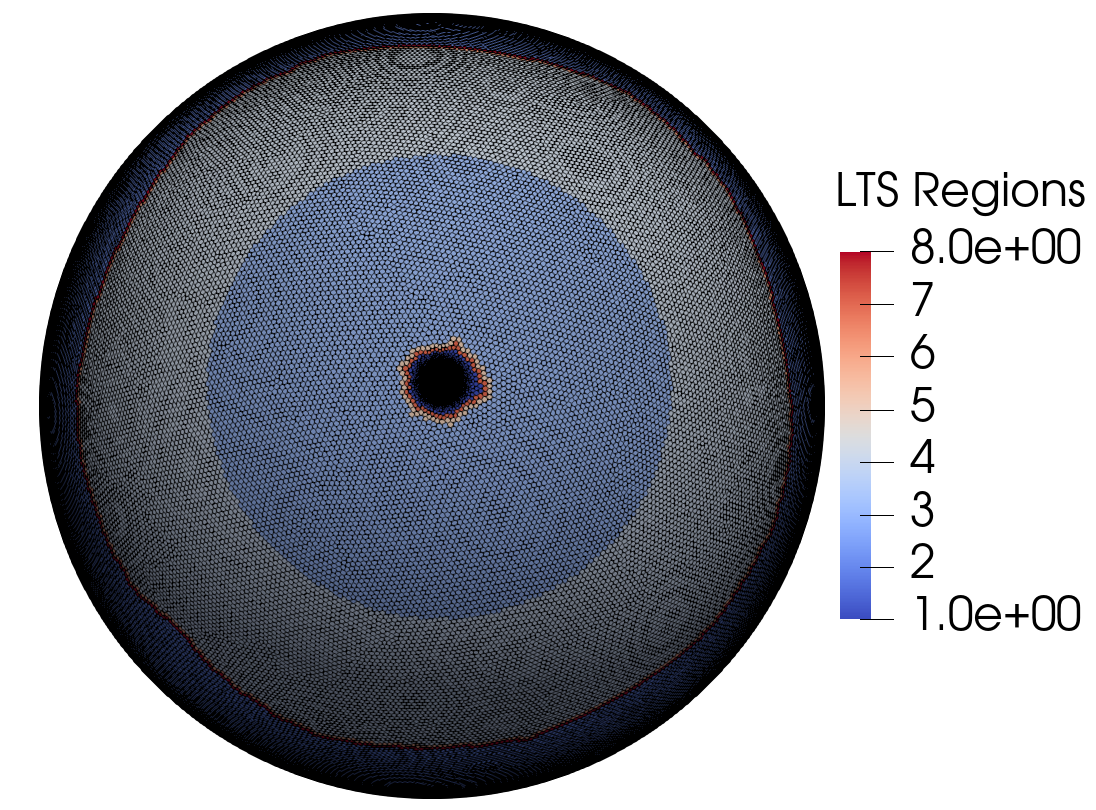}\\
      \includegraphics[scale=0.35]{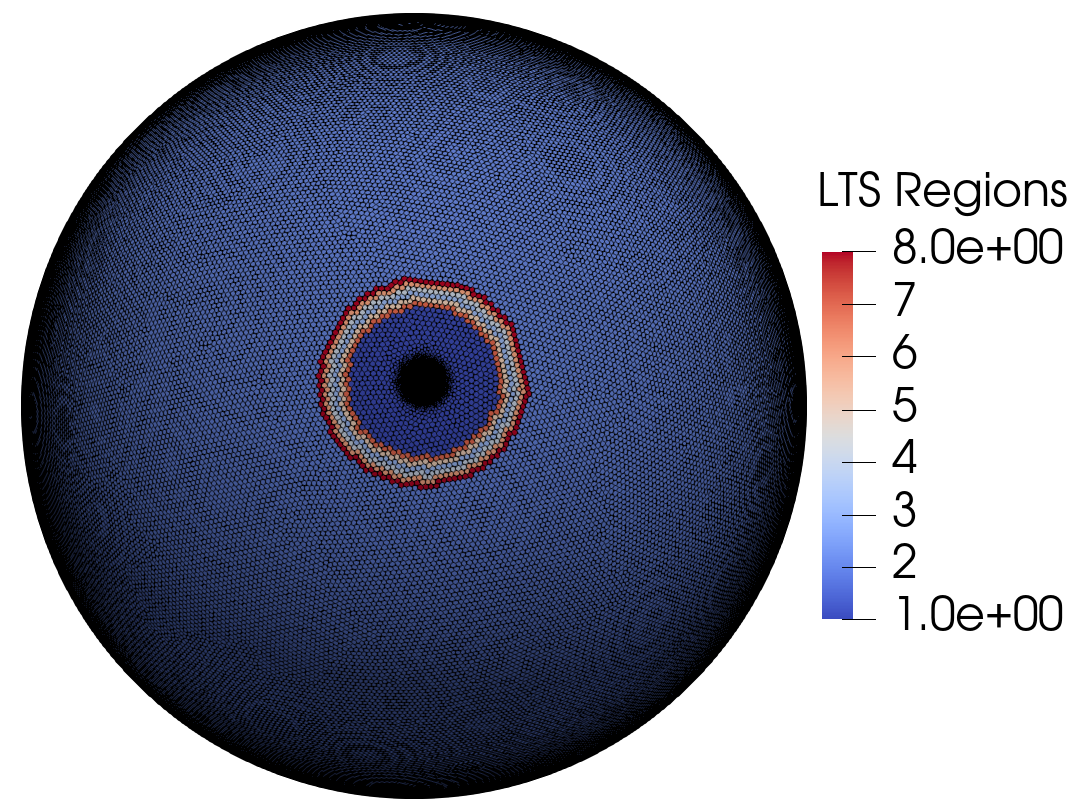}
\caption{Top: LTS regions used for case A) of the speed up test. The interface layers are composed of 18 layers of cells each. Middle: LTS regions used for case B) of the speed up test. The interface layers are composed of 38 layers of cells each. Bottom: LTS regions used for case C) of the speed up test. The interface layers are composed of 1 layer of cells each.}
   \label{fig:11}
\end{figure}

We conclude this section with what is probably the most important test to highlight the advantage of using LTS schemes on variable resolution grids.
Namely, we consider a grid with 120800 cells, where the ratio between the area of the largest cell and the area of the smallest cell is a big number, and a relatively small portion of the cells of the mesh compose the fine LTS region. 
Namely, let $A_{ls}$ and $C_{cf}$ be defined as
\begin{align}
A_{ls} = \dfrac{\mbox{area of the largest cell}}{\mbox{area of the smallest cell}}, \qquad
C_{cf} = \dfrac{\mbox{number of coarse time-step cells}}{\mbox{number of fine time-step cells}}.
\end{align}
\rev{For the mesh in consideration, the fine resolution is $2.5$ km and the coarse is $77.5$ km, and $A_{ls}=1340.273$. On the other hand, the value of $C_{cf}$ depends on the LTS regions}, for which we have three different setups:
\begin{itemize}
    \item [\bf{A})] The LTS interface layers are composed of 18 layers of cells each. In this case, we have that the coarse time-step cells (coarse + the two interface layers) are 98886, and the fine time-step cells are 21914, hence $C_{cf}=4.512$. In Figure \ref{fig:11} one can see the different LTS regions for this case, considering LTS3.
        \item [\bf{B})] The LTS interface layers are composed of 38 layers of cells each. In this case, we have that the coarse time-step cells (coarse + the two interface layers) are 98844, and the fine time-step cells are 21956, hence $C_{cf}=4.502$. In Figure \ref{fig:11} one can see the different LTS regions for this case, considering LTS3.
        \item [\bf{C})] The LTS interface layers are composed of 1 layer of cells each as in Figure \ref{fig:3}. However, all the cells in the interface layers are assigned to one processor. The procedure to achieve this setup is the following: initially, the partition of the mesh is done as explained before, i.e. each MPI rank has three memory blocks and each block is assigned a balanced number of cells coming from either the fine time-step region, the coarse time-step region or the interface layers. Then, a user selected processor is assigned all the cells belonging to the memory blocks associated with interface layers of the other processors. The result of this procedure is that a single processor has all the interface layer cells in its memory block associated with the interface layers and all the other processors do not own interface cells anymore, i.e. their block associated with interface layers has no cells.  In this case, we have that the coarse time-step cells (coarse + the two interface layers) are 98198, and the fine time-step cells are 22602, therefore $C_{cf}=4.345$. 
\end{itemize}
The reason for considering case A) and B) is that, as shown in the previous test, having a sufficient number of cells in the interface layers is important for scalability, as enough cells have to be owned by the blocks handling the interface to observe it. However, having too many cells in the interface can affect the performance of the LTS schemes in terms of CPU times.
The reason is the following: during the sub-stepping procedure, i.e. the loop $k=1,\ldots,M$ in Step 3 of our formulation of the LTS schemes, the correction terms for the interface correction are computed. These terms are computed on the interface cells, hence the more interface cells there are, the heavier this computation will be. Given that the computation is done within the sub-stepping procedure, so $M$  times per iteration, it can have a toll on the overall CPU time of the LTS scheme. In fact, we will show that case A) will go faster than B) but would not scale as good. 
Case C) is considered so that the minimum number of cells in the interface layers can be used. However, as one can imagine, such a case is the less balanced in terms of load and therefore we expect it to show a somewhat erratic scalability, which is strongly dependent on how the MPI partitioning of the mesh is carried out by METIS. Moreover, we expect C) to get worse as the number of MPI ranks is increased, due to the fact that more processors will have to be in idle while the only one that has all the interface cells is working.

\begin{figure}[!t]
   \centering
   \includegraphics[scale=0.3]{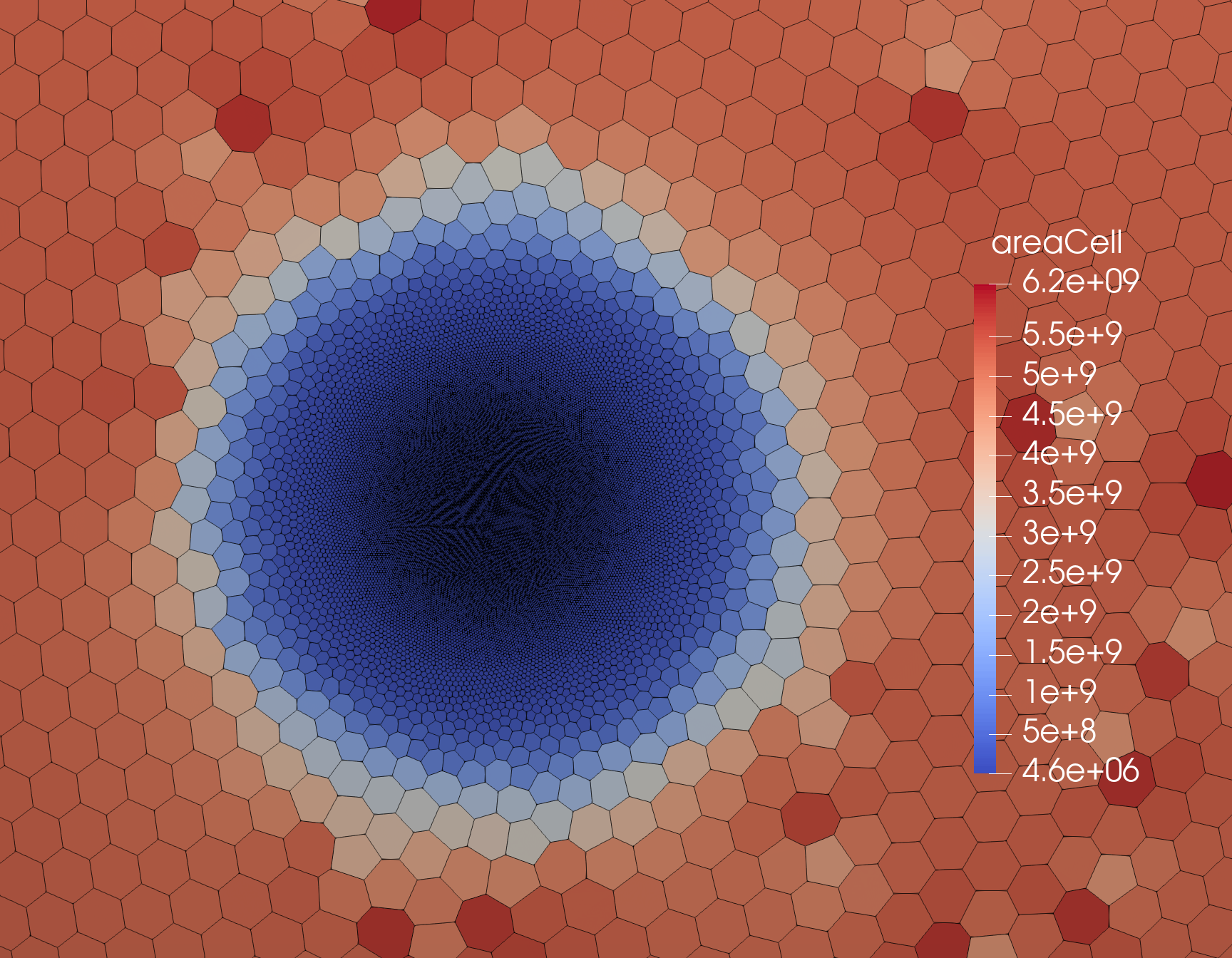}
   \includegraphics[scale=0.3]{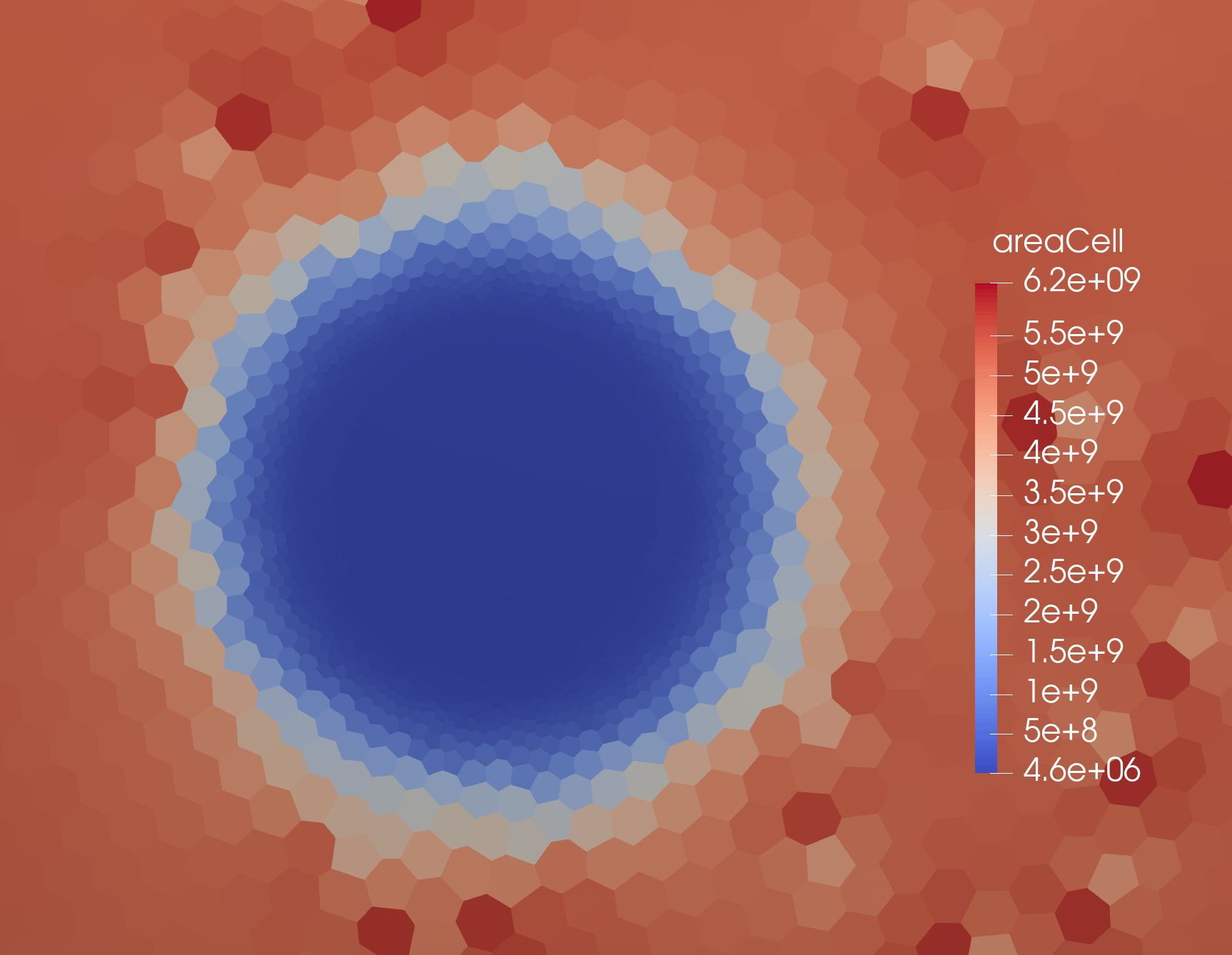}
\caption{Cell areas in proximity of the variable mesh resolution zone. Top: the mesh is displayed. Bottom: the mesh is not displayed.}
   \label{fig:10}
\end{figure}

In Figure \ref{fig:10} we are showing a zoom-in on the region of the mesh where the variable resolution occurs, and display the value of the cell areas.
Note that the mesh is the same for cases A), B) and C), so Figure \ref{fig:10} applies to all cases. The difference is that the LTS interface is placed in different places on the mesh.
%%%%%%%%%%%%
%%%%%%%%%%%%
With a mesh like the one used here, explicit time stepping algorithms without local time-stepping would find their maximum time-step constrained by the area of the smallest cell in the mesh.  With LTS, this constraint only affects the fine LTS region and much larger time steps can be selected for the coarse and interface regions, resulting in considerable computational time gains. Of course, the coarse time-step is also bounded above by a coarse CFL condition so it cannot be chosen arbitrarily. Obviously, the advantage of LTS becomes more and more evident as both $A_{ls}$ and $C_{cf}$ get bigger. 
Given that we want to explore the maximum gain in terms of CPU provided by the LTS methods, we begin by finding the largest fine time-step for the LTS schemes using $M=1$, and then progressively increase $M$ to find the largest coarse time-step.
\rev{The simulation is carried out for 1 hour and,
for the test in consideration, the maximum fine time-step for LTS3 is $\Delta t_{fine} = 9$ seconds, and the maximum coarse time-step is $\Delta t_{coarse} = 225$ seconds, which means that $M=25$, i.e. the coarse time-step is 25 times larger than the fine time-step. For LTS2, the maximum fine time-step is $\Delta t_{fine} = 2.5$ seconds, and the maximum coarse time-step is $\Delta t_{coarse} = 180$ seconds, hence $M=72$ and the coarse time-step is 72 times bigger than the fine one.
For RK4 the maximum time-step is $\Delta t = 15$ seconds.}
%%%%%%%%%%%%
%%%%%%%%%%%%
Finally, to simulate the work load that would be performed during an ocean simulation, we used 100 layers instead of 1, meaning that the computations are repeated in the same way for 100 times during each iteration. 
Results for RK4 and LTS3 are in Table \ref{table:4}, for case A) (18 layers of cells for each interface layer), B) (38 layers of cells for each interface layer), and C) (1 layer of cells for each interface layer).
\rev{Being a lower order method compared to LTS3 and RK4, LTS2 has a much lower fine time-step, and hence its performance is not a good enough candidate in terms of CPU time to beat RK4, as opposed to LTS3. For instance, with case A) and 256 MPI ranks, LTS2 shows a speed-up over RK4 of $12.6$ \% which is considerably lower than those provided by LTS3, shown in Table \ref{table:4}. Hence, to keep the focus on the very good results given by LTS3, those for LTS2 are not reported in this section.}
%%%%%%%%%%%%%%%%%%%%%%%%%%%%%%%%%%%%%%
\begin {table}[!h!]
\caption{\rev{CPU time gain for the speed up test considering case A) (18 layers of cells per interface layer), case B) (38 layers of cells per interface layer), and case C) (1 layer of cells per interface layer). The run has a duration of 1 hour. The coarse time-step for LTS3 is 225 seconds and $M=25$. The time-step for RK4 is 15 seconds. A total of 100 vertical layers is used.}}
\setlength\tabcolsep{3.25pt} % default value: 6pt
\begin{center}
\begin{tabular}{|c|c|c|c|c|c|c|c|} 
\cline{1-8}    
 & \multicolumn{7}{c|}{CPU time [s]}  \\ \hline  
   MPI ranks & RK4 & A) LTS3 ($M=25$)  &  Gain & B) LTS3 ($M=25$)  &  Gain & C) LTS3 ($M=25$)  &  Gain \\ \hline 
4 &   2132.285    &   650.889    & 69.474\% & 954.866    & 55.219\%  & 578.055 &  72.890\%\\ \hline
8 &   1265.602    &   403.021    & 68.156\% & 604.616   & 52.227\% & 343.315 & 72.873\%
\\ \hline
16 &   973.286   &  281.483    &  71.079\% & 434.624  &  55.345\%  & 272.899 & 71.961\% \\ \hline
32 &  542.968     &   150.752    &  72.235\% &225.519    &  58.465\% & 153.575& 71.715\%\\\hline
64 &   275.138    &   75.339    &  72.618\% &110.108   & 59.980\%   & 93.930 & 65.860\% \\\hline
128 &   144.191    & 43.624     &  69.746\% &57.005  & 60.465\%  & 68.600 & 52.424\% \\\hline
256 &   79.297    & 26.695     &  66.335\% &30.840  & 61.108\%  & 46.618 & 41.211\% \\\hline
\end{tabular}
\end{center}
\label{table:4}
\end{table}
%%%%%%%%%%%%%%%%%%%%%%%%%%%%%%%%%%%%%%
The advantage of using LTS3 instead of RK4 in terms of CPU time is measured as a percentage of the CPU time of RK4 and it is computed as
\begin{align}
    \mbox{Gain} = \dfrac{\Big(\mbox{CPU time of RK4} - \mbox{CPU time of LTS3}\Big) * 100}{\mbox{CPU time of RK4}}.
\end{align}
\rev{We observe that case A) is the one that shows the best performance among the three cases considered, with a percent gain that never falls below $66\%$. As expected, case B) shows gains that are slightly lower than case A), due to the increased number of computations needed by the extra cells in the interface layers. However, thanks to these cells, case B) maintains a close to linear scalability (see Table \ref{table:5}), hence no deterioration of the CPU time gain is seen in such a case going from 128 to 256. On the other hand, case A) does show a small decay of the CPU time gain going from 128 to 256 processors.}
This deterioration is expected and is caused by the fact that as the MPI ranks increase, the blocks associated with interface layer cells own a progressively smaller number of cells, hence the scalability of the LTS3 setups decreases, and so do the gains. Case C) is the case that for a small number of MPI ranks provides the best gains, but do to its erratic scalability behavior the advantage starts to fade as the number of MPI ranks grows and goes down to approximately $41\%$ for 256 MPI ranks. As mentioned above, we anticipated a deterioration in the scalability (and hence the gain) for case C) as the number of MPI ranks grows, because more processors will be in idle.
Nevertheless, case C) shows gains that are better than B) and comparable to A) for up to 32 MPI ranks after which the deterioration starts to appear.
For the reader's convenience, in Table \ref{table:5} we have reported the scalability numbers for the three LTS setups so one can clearly see how they directly influence the CPU gain in table Table \ref{table:4}.
%%%%%%%%%%%%%%%%%%%%%%%%%%%%%%%%%%%%%%
\begin {table}[!h!]
\caption{\rev{Scalability for the LTS setups in the speed up test. The run has a duration of 1 hour. The coarse time-step for LTS3 is 225 seconds and $M=25$. The time-step for RK4 is 15 seconds. A total of 100 vertical layers is used.}}
\setlength\tabcolsep{3.25pt} % default value: 6pt
\begin{center}
\begin{tabular}{|c|c|c|c|c|c|c|c|c|} 
\cline{1-9}    
 & \multicolumn{8}{c|}{CPU time [s]}  \\ \hline  
   MPI ranks &  RK4 & Ratio &A) LTS3 ($M=25$)  &  Ratio & B) LTS3 ($M=25$)  &  Ratio & C) LTS3 ($M=25$)  & Ratio \\ \hline 
4    & 2132.285 & - &    650.889    & - & 954.866    & -   & 578.055 &  - \\ \hline
8    & 1265.602 & 1.685 &    403.021    &  1.615 & 604.616   & 1.579  & 343.315 &  1.684
\\ \hline
16   & 973.286 & 1.300 &  281.483    &  1.432  & 434.624  &  1.391   & 272.899 &  1.258 \\ \hline
32  & 542.968 & 1.792   &   150.752    &  1.867 & 225.519    &  1.927  & 153.575& 1.777\\\hline
64   & 275.138 & 1.973 &   75.339    &  2.000  & 110.108   &  2.048  & 93.930&  1.635 \\\hline
128  &  144.191 & 1.908  & 43.624     &  1.727  & 57.005  &  1.931  & 68.600&  1.369 \\\hline
256  & 79.297  & 1.818    & 26.695     &  1.634  & 30.840  &  1.848  & 46.618&  1.471 \\\hline
\end{tabular}
\end{center}
\label{table:5}
\end{table}

\section{Conclusions and future work}\label{End}

In this work, we have investigated the computational aspects of the local time-stepping (LTS) schemes developed in \cite{hoang2019conservative}, focusing on implementation, parallel performance, and CPU time.
The algorithms have been implemented in the shallow water core of MPAS, with the intention of integrating the LTS procedure in the ocean core in the immediate future. Our numerical results have shown that the LTS methods implemented in MPAS achieve the expected convergence rates and possess good scalability properties, making them suitable for real world simulations of global climate models that run on hundreds of cores. 
The final test on CPU time showed that the LTS method of order 3 is capable of considerably outperforming the Runge-Kutta scheme of order 4, whereas, due to its lower order, the LTS method of order 2 can only show minor gains. Therefore, we believe that only the LTS method of order 3 should be integrated into MPAS-Ocean.
LTS algorithms present non-straightforward challenges when it comes to achieving scalability and fast parallel computations, but thanks to a careful implementation and tailored parallel partitioning, we showed that both goals can be achieved.

\section*{Acknowledgements}
The authors were supported by the Earth System Model Development (ESMD) program as part of the multi-program, collaborative Integrated Coastal Modeling (ICoM) and Energy Exascale Earth System Model (E3SM) projects, both funded by the U.S. Department of Energy, Office of Science, Office of Biological and Environmental Research.  This research used resources provided by the Los Alamos National Laboratory Institutional Computing Program, which is supported by the U.S. Department of Energy National Nuclear Security Administration under Contract No. 89233218CNA000001.

\appendix
\section{The nonlinear shallow water equations}\label{swe}
The nonlinear shallow water equations are used to describe the motion of a layer of fluid on a two-dimensional surface in a rotating reference frame, by providing conservation equations for the layer thickness $h$, and fluid velocity $\mathbf{u}$. In vector invariant form  \cite{ringler2010unified}, they read as follows:
\begin{equation}
\begin{cases}
\begin{aligned}\label{eq:swe1}
    \dfrac{\partial h}{\partial t} & = - \nabla \cdot (h\mathbf{u}),\\
    \dfrac{\partial \mathbf{u}}{\partial t} &= - \eta\mathbf{k} \times \mathbf{u}  - g \nabla(h +b) -\nabla K.
\end{aligned}
\end{cases}
\end{equation}
In the above equations, $g$ is gravity, $b$ is the bottom topography, $f$ is the Coriolis parameter, $\eta= \mathbf{k} \cdot \nabla \times \mathbf{u} + f$ is the absolute vorticity and $K=\frac{1}{2}|\mathbf{u}|^2$ is the kinetic energy.
The ratio of absolute vorticity $\eta$ and layer thickens $h$ provides an expression for the potential vorticity $q$ 
\begin{align}
    q=\dfrac{\eta}{h},
\end{align}
which can be substituted in Eq. \eqref{eq:swe1} to obtain
\begin{equation}
\begin{cases}
\begin{aligned}\label{eq:swe2}
    \dfrac{\partial h}{\partial t} &= - \nabla \cdot (h\mathbf{u}),\\
    \dfrac{\partial \mathbf{u}}{\partial t} &= - q(h\mathbf{u}^{\perp})  - g \nabla(h +b) -\nabla K.
\end{aligned}
\end{cases}
\end{equation}
The term $q(h\mathbf{u}^\perp)$, with $\mathbf{u}^\perp=\mathbf{k}\times\mathbf{u}$, is regarded as a thickness flux of potential vorticity in the direction perpendicular to the velocity field $\mathbf{u}$ \cite{ringler2010unified}. Finally, letting $\mathbf{F}=h\mathbf{u}$ define a thickness flux \cite{ringler2010unified, hoang2019conservative}, we rewrite the equations in Eq. \eqref{eq:swe2} as
\begin{equation}
\begin{cases}
\begin{aligned}\label{eq:swe3}
    \dfrac{\partial h}{\partial t} &= - \nabla \cdot \mathbf{F},\\
    \dfrac{\partial \mathbf{u}}{\partial t} &= - q(\mathbf{F}^{\perp})  - g \nabla(h +b) -\nabla K.
\end{aligned}
\end{cases}
\end{equation}
which is the form of the shallow water equations used in the remainder of the paper. Note that, similar to the definition of $\mathbf{u}^{\perp}$, $\mathbf{F}^{\perp} = \mathbf{k} \times \mathbf{F}$.
%%%%%%%%%%%%%%%%%%%%%%%%%%%%%%%%%%%%%%%%%%%%%%%%
%%%%%%%%%%%%
\begin{figure}[!t]
   \centering
   \includegraphics[scale=0.45]{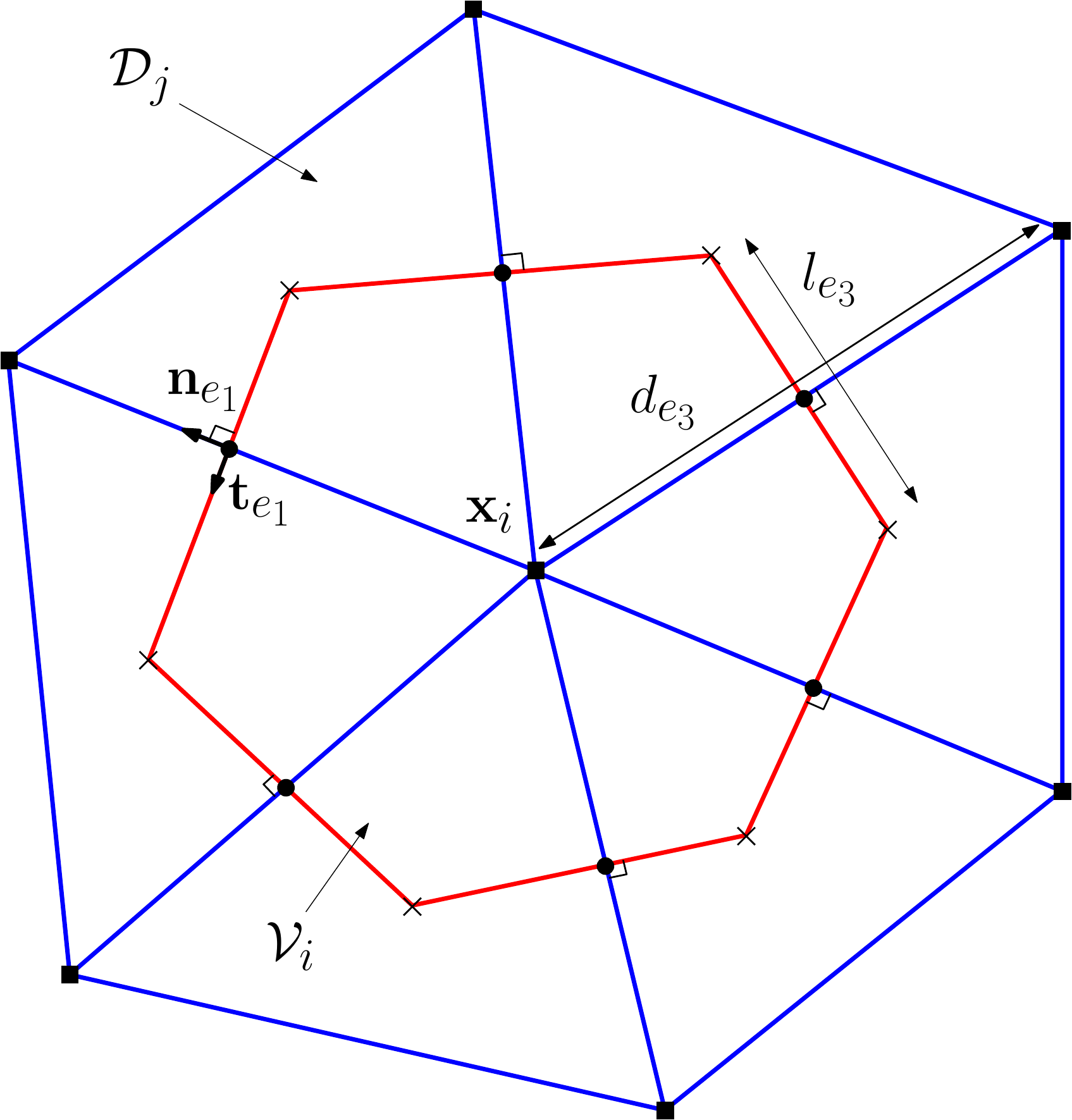}
\caption{Schematics of a Voronoi cell (red polygon) and surrounding dual Delaunay mesh (blue triangles).}
   \label{fig:1}
\end{figure}
%%%%%%%%%%%%
\section{Discretization with TRiSK}\label{TRiSK}
The nonlinear shallow water equations in Eq.\eqref{eq:swe3} are discretized on a Voronoi tesselation $\mathcal{V}=\{\mathcal{V}_i\}$ according to the TRiSK scheme, which employs a so called C-grid for the placement of the unknowns  \cite{ringler2010unified}. 
A Voronoi tesselation can be obtained by positioning a number of seed points on the surface of the sphere, then forming cells $\mathcal{V}_i$ (or regions) by grouping together all points of the surface who are closest to a specific seed point, which then identifies a cell of the tesselation \cite{ju2011voronoi}.
The special case of a SCVT is obtained after requiring that each seed point is also the centroid of its Voronoi region $\mathcal{V}_i$.
As pointed out in \cite{ringler2013multi}, SCVT grids are desirable for two major reasons: first, they allow the user to control the resolution of the mesh through the specification of desired mesh-density functions. Second, the mesh quality increases with an increasing resolution, i.e. for a fixed mesh-density function, an increasing number of grid points produces Voroni cells that tend to perfect hexagons.
The primal mesh $\mathcal{V}$ uniquely identifies a dual mesh $\mathcal{D}=\{\mathcal{D}_j\}$, which is a Delaunay triangulation. 
The relation between $\mathcal{V}$ and $\mathcal{D}$ is the following, see Figure \ref{fig:1}: the vertices of the dual triangles $\mathcal{D}_j$, marked with boxes, are interior points of the primal cells $\mathcal{V}_i$; the vertices of the primal cells, marked with crosses, are the interior points of the dual triangles. From now on, with a slight abuse of terminology, interior points of either dual or primal cells will be referred to as centers. For a given primal cell $\mathcal{V}_i$, the location of its center is denoted by $\mathbf{x}_i$. Vertices of primal cells are denoted by $\mathbf{x}_v$.
As in \cite{ringler2010unified}, we also assume that the line segments joining primal cell centers  intersect line segments joining dual triangle centers, and the segments form a right angle. For a given primal edge $e$, we denote the location on $e$ of such an intersection point with $\mathbf{x}_e$.  It is worth observing that the intersection point between the aforementioned line segments will not in general be on the primal cell edge midpoint, i.e. $\mathbf{x}_e$ will not in general be the midpoint of $e$. A modification of the TRiSK scheme where the intersection does actually occur at the midpoints has been presented in \cite{peixoto2016accuracy}, but is not considered in this work. Note that the line segments joining dual cell centers are primal cell edges $e$, whose lengths are denoted by $l_e$.
The length of the dual edge intersecting a primal edge $e$ is denoted by $d_e$, see Figure \ref{fig:1}.
At each $\mathbf{x}_e$, a local reference frame is considered, with a unit normal vector $\mathbf{n}_e$ parallel to the line segment connecting primal cell centers, and a tangential vector $\mathbf{t}_e = \mathbf{k} \times \mathbf{n}_e$, as shown in Figure \ref{fig:1}.

Within this characterization of the mesh, the fluid thickness, bottom topography and kinetic energy are located at primal mesh centers, i.e. $\mathbf{x}_i$, the points marked with a box in Figure \ref{fig:1}.
On the other hand, the normal component of the velocity $\mathbf{u}$ is specified at the primal cell edges, i.e. $\mathbf{x}_e$, the points marked with a dot in Figure \ref{fig:1}. This C-grid characterization of the velocity field that does not specify the tangential component requires such a component to be recovered when needed, unlike for instance for B-grids, where both components of the velocity field are specified at the primal cell vertices \cite{bouillon2009elastic}.
The potential vorticity $q$ is instead defined at the primal cell vertices, i.e. $\mathbf{x}_v$, the points marked with a cross in Figure \ref{fig:1}.
In order to discretize the continuous operators in Eq. \eqref{eq:swe3} on the C-grid, it is necessary to introduce finite sets of cells and edges as follows \cite{ringler2010unified,hoang2019conservative}:
\begin{equation}
\begin{aligned}
    \mbox{For any primal cell $\mathcal{V}_i$}&: \mathcal{EC}(i) = \{\mbox{set of primal edges $e$ of $\mathcal{V}_i$.}\}\\
    \mbox{For any primal cell vertex $v$}&:
    \mathcal{EV}(v)=\{\mbox{set of primal edges $e$ sharing vertex $v$.}\}\\
        \mbox{For any primal cell edge $e$}&:
    \begin{cases}
    \mathcal{CE}(e)=\{\mbox{the two primal cells $\mathcal{V}_i$ and $\mathcal{V}_{i'}$ that own $e$.} \}\\
    \mathcal{VE}(e)=\{\mbox{the two vertices $v$ and $v'$ that are the endpoints of $e$.} \}\\
    \mathcal{EE}(e)=\{\mbox{set of primal edges $e'\neq e$ such that $e'$ is owned by cells in $\mathcal{CE}(e)$.}\}
    \end{cases}
 \end{aligned}   
\end{equation}
\begin{remark}
The choice of names for the above sets might seem a bit obscure for a reader not familiar with the MPAS framework. The truth is that such names are shorthand notations for MPAS functions that return the elements of these sets. For instance, $\mathcal{EV}$ stands for "edges on vertex", which is the MPAS function that returns the primal edges sharing a primal input vertex. 
\end{remark}
Considering the construction presented in this section, the equations in Eq. \eqref{eq:swe3} are written in discrete form as follows
\rev{
\begin{equation}
\begin{cases}
\begin{aligned}\label{eq:sweDiscrete}
    \dfrac{\partial h_i}{\partial t} &= - \Big[\nabla \cdot F_e\Big]_i,\\
    \dfrac{\partial u_e}{\partial t} &= - F_e^{\perp}[q]_{v \rightarrow e}  - \Big[g \nabla(h_i +b_i) + \nabla K_i\Big]_e,
\end{aligned}
\end{cases}
\end{equation}
}
note that the equation for the layer thickness is discretized at the primal cell centers, whereas the discretized momentum equation is obtained for the primal edge locations, and is derived by taking the dot product between $\mathbf{n}_e$ and the momentum equation in Eq. \eqref{eq:swe3} at the edge locations, for all primal edges $e$.
The notation $[\cdot]_{v \rightarrow e}$ refers to an operator that provides at an edge location the value of a quantity that would be defined at a vertex.
For the potential vorticity, such an operator is defined as \cite{ringler2010unified}
\begin{equation}\label{eq:discPV}
[q]_{v \rightarrow e} = \dfrac{1}{F_e^{\perp}} \Big( \sum\limits_{e'\in \mathcal{EE}(e)} w_{e,e'}\dfrac{l_{e'}}{d_e}F_{e'} (\dfrac{\widetilde{q}_e + \widetilde{q}_{e'}}{2})\Big), \qquad  \widetilde{q}_e = \dfrac{1}{2}  \sum\limits_{v\in \mathcal{VE}(e)} q_v.
\end{equation}
The expression of $\widetilde{q}_e$ guarantees energy conservation \cite{hoang2019conservative}, alternative formulations are presented in \cite{ringler2010unified}.
A similar operator is used within the discrete thickness flux $F_e =  [h]_{i\rightarrow e} u_e$, i.e. the operator  $[\cdot]_{i\rightarrow e}$ maps the layer thickness from a center to an edge.
 As remarked in \cite{hoang2019conservative}, the specific choice of operator depends on the physics of the problem, in this work we consider 
 \begin{equation}
     [h]_{i \rightarrow e} = \dfrac{1}{2}\sum\limits_{i \in \mathcal{CE}(e)}h_i.
 \end{equation}
The term $F_e^{\perp}$ is described by the following expression \cite{thuburn2009numerical}
\begin{equation}\label{eq:Feperp}
    F_e^{\perp}= \sum\limits_{e'\in \mathcal{EE}(e)} w_{e,e'}\dfrac{l_{e'}}{d_e}F_{e'}.
\end{equation}
For two primal edges $e$ and $e'$, let $\mathcal{V}_{i^*}$ be the primal cell that owns both $e$ and $e'$, i.e. $\mathcal{V}_{i^*} = \mathcal{CE}(e) \cap \mathcal{CE}(e')$. Let $V_{i^*}$ be the set of vertices ordered in either clockwise or counterclockwise direction found when moving from edge $e'$ to edge $e$ along $\mathcal{V}_{i^*}$. Let $v^*$ be the last of such vertices. Note that $v^*$ will be either one of the vertices of $e$.
Then, for any $e' \in \mathcal{EE}(e)$, the weights in Eq. \eqref{eq:discPV} and \eqref{eq:Feperp} are defined as \cite{thuburn2009numerical}: 
\begin{equation}\label{eq:weights}
w_{e,e'}=
     n_{e',i^{*}}t_{e,v^*}\Big(\sum\limits_{v \in V_{i^*}}\dfrac{A_{i^{*},v}}{A_{i^{*}}} -\dfrac{1}{2}\Big). \end{equation}
A detailed description of the symbols employed in Eq. \eqref{eq:weights} is in order:
the terms $n_{e,i}$ and $t_{e,v}$ are sign corrections defined as follows:
\begin{equation}
n_{e,i}=
    \begin{cases}
    1, \qquad \mbox{if $\mathbf{n}_e$ is an outward normal for $\mathcal{V}_i$}\\
    -1, \qquad \mbox{otherwise}.
    \end{cases}
\end{equation}
\begin{equation}
t_{e,v}=
    \begin{cases}
    1, \qquad \mbox{if $v$ is in the direction $\mathbf{k} \times \mathbf{n}_e$ (i.e. $v$ is at the left end of edge $e$ \cite{thuburn2009numerical}).}\\
    -1, \qquad \mbox{otherwise}.
    \end{cases}
\end{equation}
\begin{remark}
The direction of the normal vectors $\mathbf{n}_e$ is arbitrary, hence in MPAS, the value of $n_{e,i}$ is decided according to the following convention: if $\mathcal{V}_{i_1}$ and $\mathcal{V}_{i_2}$ are the primal cells in $\mathcal{CE}(e)$, and if $i_1 > i_2$, then we have $n_{e,i_1}=1$ and $n_{e,i_2}=-1$.
\end{remark}
The quantity $A_{i^*}$ represents the area of cell $\mathcal{V}_{i^*}$, whereas the quantities $A_{i^*,v}$ are informally referred to as the {\it kite areas}, and are the areas of polygonal subsets of $\mathcal{V}_{i^*}$ whose vertices are the cell center $\mathbf{x}_{i^*}$, the vertex $\mathbf{x}_v$, and the edges points $\mathbf{x}_{e_1}$ and $\mathbf{x}_{e_2}$ such that $v \in \mathcal{VE}(e_1)\cap \mathcal{VE}(e_2)$, please see \cite{thuburn2009numerical} for further details.
The discrete kinetic energy $K_i$ is computed as follows \cite{ringler2010unified,hoang2019conservative}:
\begin{equation}
    K_i = \dfrac{1}{4 \, A_i} \sum\limits_{ e \in \mathcal{EC}(i)} l_e d_e u_e^2,
\end{equation}
where once again $A_i$ represents the area of the primal cell indexed by $i$.
Finally, it is left to discuss the discretization of the operators in Eq. \eqref{eq:sweDiscrete}:
the divergence term $[\nabla \cdot F_e]_i$ is given by
\begin{align}
    \Big[\nabla \cdot F_e\Big]_i = \dfrac{1}{A_i}\sum\limits_{ e \in \mathcal{EC}(i)} n_{e,i}l_eF_e,
\end{align}
and the gradient terms of the layer thickness in Eq. \eqref{eq:sweDiscrete} is expressed as
\begin{equation}
    \Big[\nabla h_i\Big]_e = \dfrac{1}{d_e}\sum\limits_{ i \in \mathcal{CE}(e)} - n_{e,i} h_i,
\end{equation}
and similarly for the bottom topography $b$ and kinetic energy $K$.
%%%%%%%%%%%%%%%%%%%%%%%%%%%%%%%%%%%%%%%%%%%%%%%%

\section*{References}
\bibliographystyle{elsarticle-num}
\bibliography{bib.bib}

\end{document}